\documentclass[final,onefignum,onetabnum]{siamonline190516}
\usepackage{amsmath}
\usepackage{pinlabel}
\usepackage{mathtools} 
\usepackage[T1]{fontenc}

\usepackage{physics}
\usepackage{multirow}
\usepackage{tikz}
\tikzset{>=latex}
\usetikzlibrary{positioning}
\usepackage{float}

\usepackage{lipsum}
\usepackage{amsfonts}
\usepackage{graphicx}
\usepackage{epstopdf}
\usepackage{algorithmic}
\ifpdf
  \DeclareGraphicsExtensions{.eps,.pdf,.png,.jpg}
\else
  \DeclareGraphicsExtensions{.eps}
\fi

\newsiamremark{remark}{Remark}
\newsiamremark{hypothesis}{Hypothesis}
\newsiamremark{question}{Question}
\crefname{hypothesis}{Hypothesis}{Hypotheses}
\newsiamthm{claim}{Claim}

\title{Wallpaper Groups and Auxetic Metamaterials
}

\author{Brendan Burns Healy\thanks{Department of Mathematical Sciences, University of Wisconsin--Milwaukee, Milwaukee, WI 
  (\email{healyb@uwm.edu}).}
  \and Aparna Deshmukh\thanks{Department of Civil Engineering, University of Wisconsin--Milwaukee, Milwaukee, WI.}
\and Elliott Fairchild$^\dagger$
\and Caroline J. Merighi\thanks{Department of Mathematical and Statistical Sciences, Marquette University, Milwaukee, WI}.
\and Konstantin Sobolev$^\ddagger$}

\headers{Wallpaper Groups and Auxetic Metamaterials}{B. B. Healy, A. Deshmukh, E. Fairchild, C. J. Merighi, and K. Sobolev}

\usepackage{amsopn}

\newsiamremark{convention}{Convention} 

\usepackage{multicol}
\usepackage{changepage}

\begin{document}

\maketitle

\begin{abstract}
We examine a fundamental material property called Poisson's ratio, which establishes the relationship for the relative deformation of a physical system in orthogonal directions. Architects and engineers have designed advanced systems using repeating patterns that can potentially exhibit auxetic behavior, which is the property of having a negative Poisson's ratio. Because two-dimensional cross sections of each of these patterns has an associated wallpaper group, we can look for useful correlations between this geometric information and the two-dimensional response as defined by Poisson's ratio. By analyzing the data, we find two properties of the wallpaper group that correlate with more effective Poisson's ratio required for applications. This paper also contains an introduction to wallpaper groups and orbifold notation and an appendix contains some literature references recreated to use our preferred notation.
\end{abstract}



\section{Introduction}
It has long been an area of interest in materials science and engineering to study how particular small-scale arrangements of a material can contribute to aggregate properties seemingly in contradiction with the large-scale properties of the bulk material itself. Such arrangements of material are called \emph{metamaterials}, and many classes of these metamaterials have proven useful in countless applications. 

While much of the scientific literature focuses on metamaterials whose target properties are photonic and electromagnetic in nature, this paper will concentrate on mechanical properties of these systems. Metamaterials designed for mechanical optimization are useful for their resistance to permanent deformation and durability under strain. Mechanical metamaterials are of interest for advanced design and construction of concrete/reinforced concrete structures such as buildings and roads as well as medical applications like durable and flexible prosthetics. For more examples, see \cite[Section~1.3]{applications}.



One class of particular interest is those systems which undergo a particular type of deformation, or shape change, when force is applied. In two dimensions, a system under stress, or force in a particular direction, will usually see points shift (undergo \emph{strain}) in that direction, as well as the orthogonal (\emph{transverse}) direction. 

\begin{definition}[2-dimensional Poisson's ratio]
Let $x$ be the direction in which the force is applied to a metamaterial system, and $z$ be an orthogonal direction. Then we let $\epsilon_x, \epsilon_z$ denote how much deformation occurs in each of those directions. \footnote{More precisely, given a displacement field $\textbf{u}$, we have $\epsilon_x= \frac{\partial \textbf{u}_x}{\partial x}$ and $\epsilon_z= \frac{\partial \textbf{u}_z}{\partial z}$. }
The \emph{Poisson's ratio} of the material in these directions is given by $\nu_{xz} := - \frac{\epsilon_{z}}{\epsilon_{x}}$. A material is \emph{auxetic} in these directions if $\nu_{xz} < 0$. 
\end{definition}
\begin{figure}[h]
\begin{center}
\includegraphics[scale=0.5]{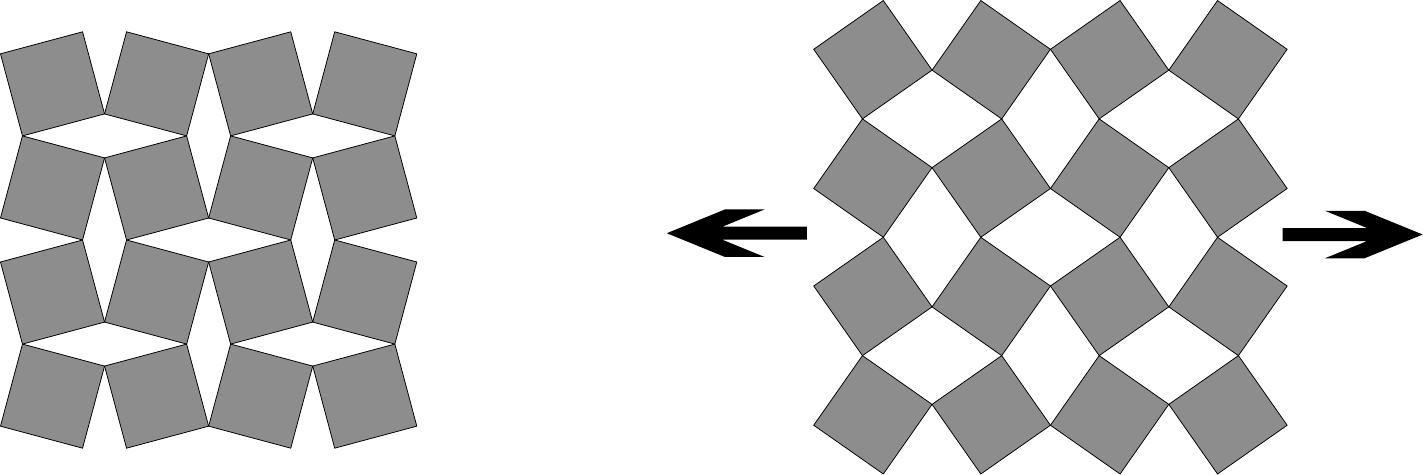}
\end{center}
\caption{An example of an auxetic system in the relaxed state (left) and under strain in the horizontal direction (right).}
\label{fig:ex}
\end{figure}

We note here that it is possible for a material to be auxetic in some directions, but not auxetic in other directions. In other words, both \emph{isotropic} auxetic metamaterials, those for whom the Poisson's ratio in 2-dimensional space is independent of the directions $x$ and $z$, and \emph{anisotropic} auxetic metamaterials, those for whom the ratio depends on the directions $x,z$, exist. For more background on auxetic metamaterials, the reader is referred to \cite{review1}, \cite{review}, \cite{lim}.

As most auxetic systems are periodic in nature, we can gain insight into them by studying their inherent periodic geometries. This idea has been pursued by Borcea and Streinu in their development of an expansive mathematical theory of auxetic behavior in arbitrary dimensions by studying one-parameter deformations of periodic lattices (see  \cite{borcea2015geometric}, as well as \cite{streinu2}, \cite{streinu3}). Here we will zoom in on a particular property of 2-dimensional auxetic systems; we study the associated symmetry groups, called \emph{wallpaper groups}, which will be introduced and discussed in the next section. Some consideration of the wallpaper group for an auxetic system has already been studied in \cite{wallpaper_first}. In this paper, we extend this relationship by analyzing how the wallpaper group corresponding to an auxetic system contributes to its deformation properties. 


Section~\ref{sec:Background} is a self-contained introduction to wallpaper groups, including orbifold notation. We discuss the geometry of the topological space obtained by taking a quotient of the Euclidean plane by a given wallpaper group considered as a subgroup of the isometries of the plane.

In Section~\ref{sec:Hypothesis} we outline one possible metric for wallpaper groups as a way of correlating the symmetry of a metamaterial system and its auxetic properties, which themselves will be primarily be measured using Poisson's ratio. In order to do this, we discuss how one wallpaper group may be considered to represent more symmetry than another. We also discuss how motivation for the guiding questions of the present article are highly influenced by findings of Grima, Mizzi, Azzopardi, and Gatt, who in \cite{random} demonstrate that a purposeful loosening of the symmetries of a system negatively impact its auxetic performance.

In Section~\ref{sec:data} we list the data we are collecting from the literature to examine for correlations. For each paper we draw from, we give a quick summary of the relevant portions of the paper, in order to contextualize the systems under analysis. Finally, Section~\ref{sec:conclusions} summarizes the conclusions we draw from this data, and Appendix~\ref{sec:Further} lists a few further questions of interest to the authors. In particular we find that reflective symmetry has a beneficial effect on Poisson's ratio, and that rotational symmetry of some orders are better than others.

\section{Wallpaper Groups}
\label{sec:Background}

Any two-dimensional pattern that can be extended infinitely to fill the Euclidean plane can be classified based on the types of symmetry present in that pattern. We can describe which type(s) of symmetry are present in the pattern by identifying the \emph{symmetry group}, or group of isometries, associated to that pattern. Recall that an isometry is a distance preserving homeomorphism of the plane $\mathbb{E}^2$.

All isometries of the Euclidean plane are of one of the following four types: translation, rotation, reflection, or glide reflection. A \emph{wallpaper group} is a symmetry group of the plane that contains two (linearly) independent translations. Wallpaper groups are classified by how many (if any) of each type of the other three symmetries are present in the pattern

There are exactly 17 wallpaper groups up to abstract isomorphism, a fact proved by Evgraf Fedorov in 1891 \cite{Fedorov} and independently by George P\'{o}lya in 1924 \cite{polya}. 
One notable part of this classification is that a wallpaper group can only contain rotations by $60^\circ, 90^\circ, 120^\circ,$ or $180^\circ$. This comes from the crystallographic restriction theorem (see, for instance, \cite{coxeterbook}), which in two dimensions can be formulated as: Every discrete isometry group of the Euclidean plane with translations spanning the whole plane, any isometry of finite order\footnote{Recall that the order of a group element $g$ is the smallest possible integer such that $g^n$ is the identity (trivial) element. Note that any isometries of order 1 are trivial so we do not need to keep track of them.} must have order 1, 2, 3, 4 or 6. Note that all reflections have order $2$. Rotations of order 2 are $180^\circ$ rotations, order $3$ are $120^\circ$ rotations, order $4$ are $90^\circ$ rotations, and order $6$ are $60^\circ$ rotations.

\subsection{Fundamental Domain} The fundamental domain of a wallpaper pattern is the smallest piece of the pattern that will create the entire pattern when moved around by the isometries of the wallpaper group. You can think of it as the smallest ``tile'' that you would need to cover the whole plane without overlaps of gaps. The ways you would need to move the tile around by translations, rotations, reflections, and/or glide reflections are the elements of the wallpaper group associated with the pattern. 

\begin{figure}[h]
\begin{center}
\includegraphics[scale=0.5]{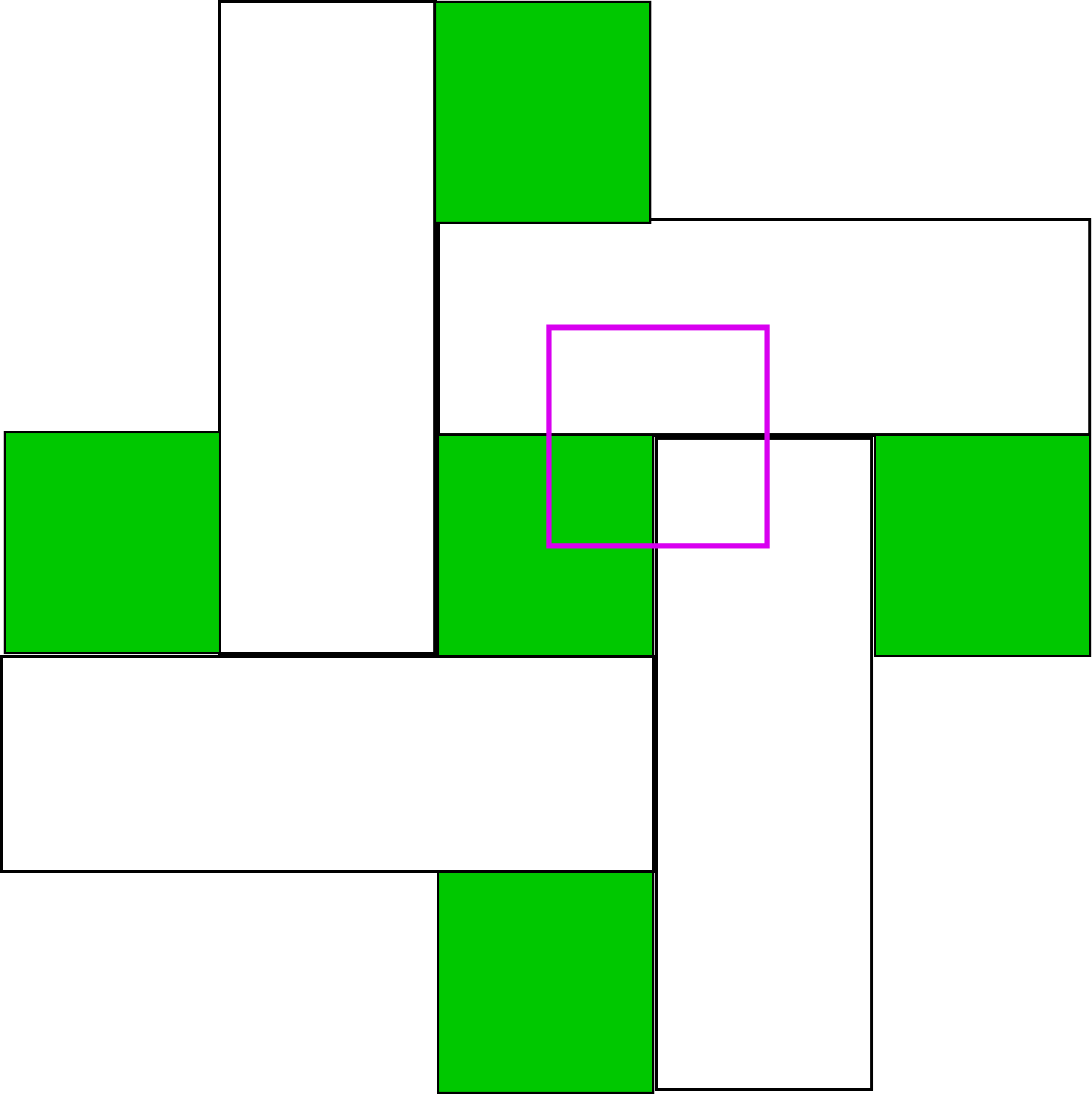}
\end{center}
\caption{The region outlined in purple is a fundamental domain for this pattern.}
\label{fig:funddom}
\end{figure}

The wallpaper group acts on the fundamental domain by identifying its vertices and edges. If we quotient $\mathbb{E}^2$ by these identifications, we get a particular \emph{orbifold}, which is a topological object encoding the wallpaper group data. Essentially, one can take the fundamental domain and ``fold'' and ``glue'' it based on the symmetries of the group. The shape that results is our orbifold. For wallpaper groups, there are 17 possible (2-)orbifolds that can result from this process\footnote{When starting with the plane, that is. We caution the reader that there are many other two-dimensional orbifolds, but their universal covers are not $\mathbb{E}^2$.}, and they correspond exactly to the 17 wallpaper groups. We can characterize which quotient space (orbifold) we get using \emph{orbifold notation}\footnote{Other notation systems for wallpaper groups are also in use, including Coxeter and crystallographic notation.}, which in turn gives notation for the wallpaper groups.

\begin{figure}[h]
\begin{center}
\includegraphics[scale=0.6]{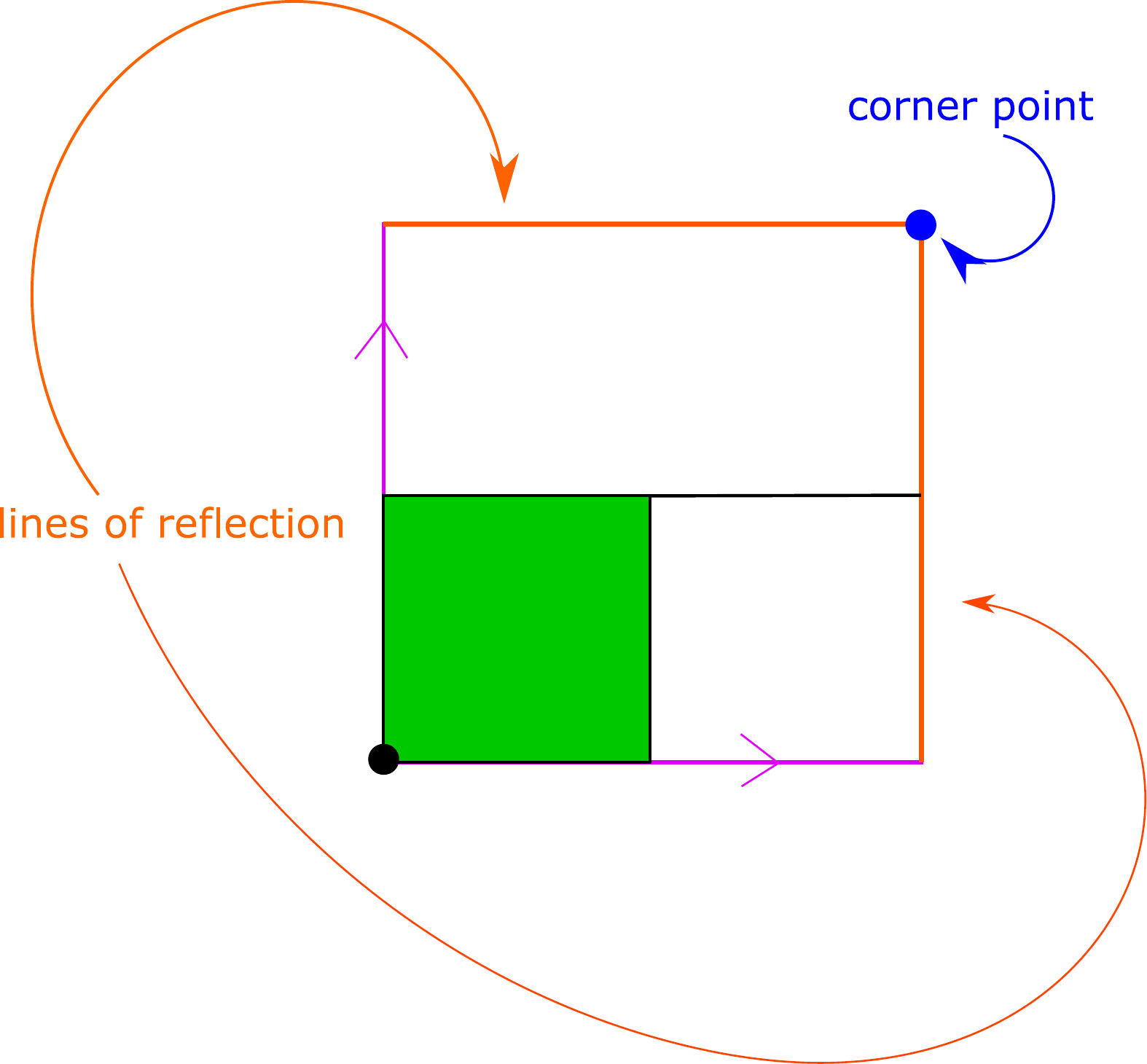}
\end{center}
\caption{Orbifold identifications in the fundamental domain of the pattern shown in Figure~\ref{fig:funddom}.}
\label{fig:orbifoldcover}
\end{figure}

\begin{figure}[h]
\begin{center}
\includegraphics[scale=0.6]{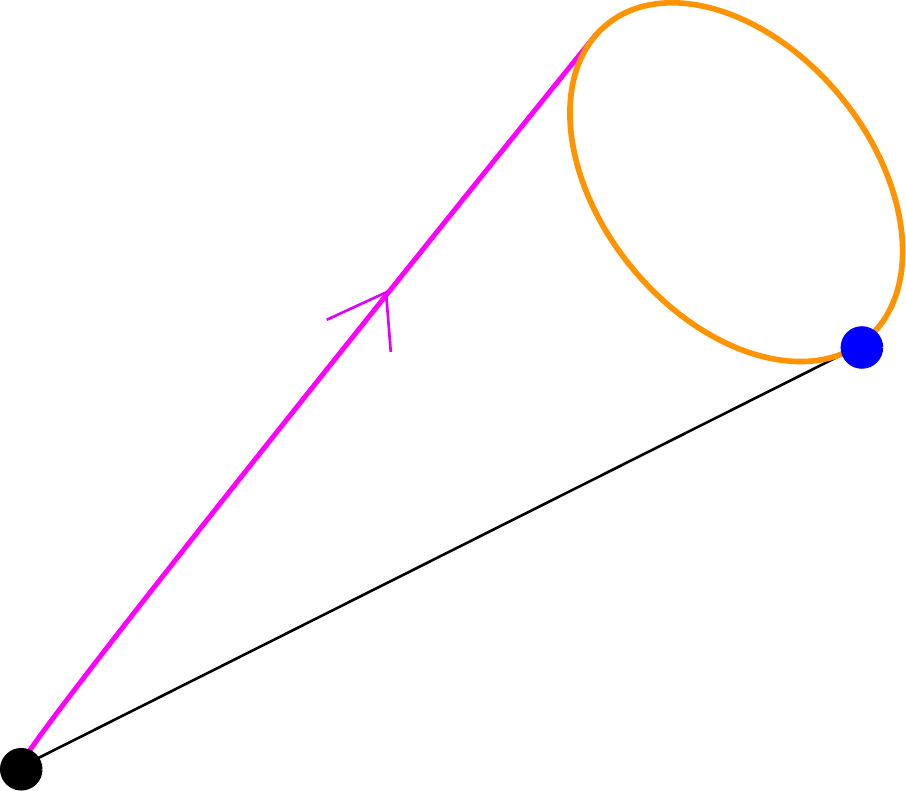}
\end{center}
\caption{The quotient orbifold $4^* 2$ (green region left uncolored).}
\label{fig:orbifoldquotient}
\end{figure}

For example, the symmetries of the pattern shown in Figure \ref{fig:funddom} are an order 4 rotation and two orthogonal lines of reflection (mirror symmetry). The quotient of this space by the wallpaper group that consists of these symmetries (aka the orbifold) has a cone point of order $4$ (this comes from the order $4$ rotation) and a corner point of order $2$ (which comes from the intersection of the two lines of reflection. Visual representations of this orbifold are shown in Figures \ref{fig:orbifoldcover} and \ref{fig:orbifoldquotient}.

\subsection{Orbifold Notation}
There are several different conventions for naming the 17 wallpaper groups. In this paper, we use orbifold notation \cite{orbifold}. It is a natural choice because it makes it relatively easy to ``see'' the topology. In this notation system, each wallpaper group has a unique \emph{orbifold symbol} that consists of a string of characters. Each character (and its location in the string) indicates a characteristic of the orbifold. 

While orbifold notation can be used to describe more than just wallpaper groups, here we will only discuss the characters used in the orbifold symbols of the 17 wallpaper groups. These characters are:
\begin{itemize}
    \item o : denotes the group with only translational symmetry
    \item $\times$ : denotes a glide reflection (a cross-cap on the orbifold)
    \item $^*$ : indicates whether numerals denote cone points or corner points
    \item $2, 3, 4, 6$ : denote the order of each cone/corner point
\end{itemize}

 Each numeral after an asterisk corresponds to a \emph{corner point}, which is a point of rotation that lies on the boundary of the quotient orbifold (see the blue point in Figure~\ref{fig:orbifoldquotient}. Each numeral before an asterisk (or in an orbifold symbol that does not contain an asterisk) corresponds to a \emph{cone point}, which is a point of rotation that does not lie on a line of mirror symmetry, i.e. is in the interior of the orbifold (see the black point in Figure~\ref{fig:orbifoldquotient}. The numeral indicates the order of the rotation. For example, the orbifold symbol for the orbifold in Figure \ref{fig:orbifoldcover}. is $4^*2$. The $4$ before the $^*$ tells us there is a cone point of order 4 (a center of rotation by $90^\circ$ that does not lie on a line of reflection). The $2$ after the $^*$ tells us there is one corner point of order $2$ (a center of rotation by $180^\circ$ that does lie on a line of reflection).

 \section{Symmetry and Auxetic Properties}
 \label{sec:Hypothesis}
 
The goal of this paper is to explore whether useful correlations exist between the geometric data of a metamaterial system and its auxetic properties. This investigation is conducted primarily via re-analyzing data on auxetic systems found in the mechanical engineering literature. We summarize auxetic systems and Poisson's ratios that have been investigated in prior papers and then identify the wallpaper group associated to each system. We then look for relationships between the geometric data of the wallpaper group and the Poisson's ratio of the system, and discuss our conclusions in Section~\ref{sec:conclusions}. Through this analysis we provide some preliminary insights into the geometric properties correlated with better performing auxetic systems, in order to aid in the design of effective metamaterials in the future.
 
To motivate this examination, we turn to one existing result on how the Poisson's ratio of an auxetic system varies with its regularity. In \cite{random}, Grima, Mizzi, Azzopardi, and Ruben find that the auxetic properties of metamaterial system degrade with a controlled destruction of the symmetry. In particular, they begin with an ordered, highly symmetric repeating two-dimensional pattern of cuts, then allow for those cuts to rotate randomly. As this random variation is allowed to increase, they find that the Poisson's ratio monotonically increases (which is a less desirable behavior). While the author's conclusion is that minor variations have little effect on the system's performance, we may take the monotonic behavior as evidence for our hypothesis that more symmetry yields better systems.

Before we can explicitly state our hypothesis, we must carefully state what we mean by ``more symmetric''. Intuitively, a wallpaper group which is the same as another except that it contains reflections, for example, should be considered ``more symmetric''. Similarly, a pattern is ``less symmetric'' than another if it only has 3-fold, rather than 6-fold, rotation. To do this comparison formally, we refer to work of Coxeter \cite{generators}. In this paper, Coxeter investigates (among many other things) the relationships that the wallpaper groups have to each other in an abstract algebraic sense. He answers the question as to which groups abstractly embed as subgroups of the other. This relation almost forms a poset (Partially Ordered SET) structure on the wallpaper groups, except that unfortunately some pairs of groups mutually embed in each other. However, if we treat those groups as equivalent, we can form a partial order on the equivalence classes of these groups.

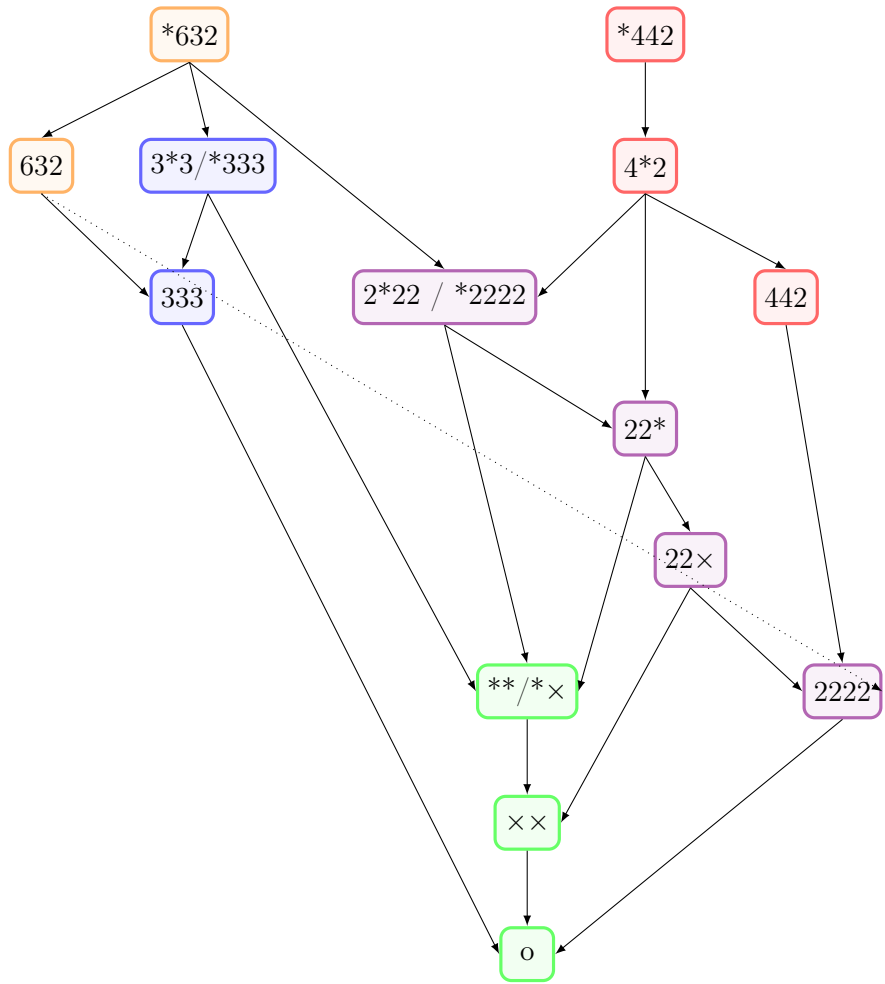
\begin{figure}[h!]
    \centering
\begin{tikzpicture}[
greennode/.style={rectangle, rounded corners, draw=green!60, fill=green!5, very thick, minimum size=7mm},
orangenode/.style={rectangle, rounded corners, draw=orange!60, fill=orange!5, very thick, minimum size=7mm},
rednode/.style={rectangle, rounded corners, draw=red!60, fill=red!5, very thick, minimum size=7mm},
bluenode/.style={rectangle, rounded corners, draw=blue!60, fill=blue!5, very thick, minimum size=7mm},
violetnode/.style={rectangle, rounded corners, draw=violet!60, fill=violet!5, very thick, minimum size=7mm},
]
\node[orangenode]   (*632)                                                  {*632};
\node[orangenode]   (632)           [below left=of *632]                    {632};
\node[bluenode]     (3*3)           [xshift=-2.2cm, below right=of *632]    {3*3/*333};
\node[bluenode]     (333)           [below right=of 632]                    {333};
\node[violetnode]   (2*22)          [below right=of 3*3]                    {2*22 / *2222};
\node[rednode]      (4*2)           [above right=of 2*22]                   {4*2};
\node[rednode]      (*442)          [above=of 4*2]                          {*442};
\node[rednode]      (442)           [below right=of 4*2]                    {442};
\node[violetnode]   (22*)           [below right=of 2*22]                   {22*};
\node[violetnode]   (22x)           [xshift=.6cm, below=of 22*]             {22$\times$};
\node[violetnode]   (2222)          [below right=of 22x]                    {2222};
\node[greennode]    (**)            [below left=of 22x]                     {**/*$\times$};
\node[greennode]    (xx)            [below=of **]                           {$\times\times$};
\node[greennode]    (o)             [below=of xx]                           {o};
\draw[->] (*632.south) -- (632.north);
\draw[->] (*632.south) -- (2*22.north);
\draw[->] (*632.south) -- (3*3.north);
\draw[->] (632.south) -- (333.west);
\draw[->] (3*3.south) -- (333.north);
\draw[->] (*442.south) -- (4*2.north);
\draw[->] (4*2.south) -- (442.north);
\draw[->] (4*2.south) -- (22*.north);
\draw[->] (4*2.south) -- (2*22.east);
\draw[->] (442.south) -- (2222.north);
\draw[->] (22*.south) -- (22x.north);
\draw[->] (22*.south) -- (**.east);
\draw[->] (3*3.south) -- (**.west);
\draw[->] (2*22.south) -- (22*.west);
\draw[->] (2*22.south) -- (**.north);
\draw[->] (22x.south) -- (2222.west);
\draw[->, dotted] (632.south) -- (2222.east);  
\draw[->] (2222.south) -- (o.east);
\draw[->] (**.south) -- (xx.north);
\draw[->] (xx.south) -- (o.north);
\draw[->] (333.south) -- (o.west);
\draw[->] (22x.south) -- (xx.east);

\end{tikzpicture}
    \caption{A partial order of (equivalence classes of) the seventeen wallpaper groups, color-coded by the highest order rotation in the group. 
    We say a group $G$ is ``more symmetric'' than another group $H$ if there is an directed path from $G$ to $H$. For example, $2^*22$ is more symmetric than $22\times$ because there is a path from $2^* 22$ to $22\times$ (via $22^*$).}
    \label{fig:poset}
\end{figure}

\begin{definition}
Let $G,H$ be wallpaper groups. We say that $G$ is \emph{more symmetric} than $H$, written $H \leq G$, if there is a path (downward) from $[G]$ to $[H]$ in Figure~\ref{fig:poset}.
\end{definition}

The use of the notation $[G]$ is to recall that we want to consider equivalent any two groups which share a box in Figure~\ref{fig:poset}. The following is a routine exercise in set theory (indeed, the proof is exactly the construction of the above figure from the information of \cite{generators}).

\begin{proposition}
Let $\mathcal{W}$ be the set of equivalence classes of wallpaper groups. Then $(\mathcal{W}, \leq)$ is a poset.
\end{proposition}

We can now formally state the primary question of interest for the present article.

\begin{hypothesis}
\label{hypothesis}
Auxetic systems associated to more symmetric wallpaper groups have lower Poisson's ratio.
\end{hypothesis}

In the following sections we will collect and analyze data from the literature to investigate Hypothesis~\ref{hypothesis}. We make the following (standard) definitions relevant to the geometry of auxetic systems.
    
   \begin{definition}
    A metamaterial system is called \emph{anisotropic} if there exist directions $x,z$ such that $\nu_{xz} \neq \nu_{zx}$. A system which is not anisotropic is said to be \emph{isotropic}.
    \end{definition}
    
       \begin{definition}
    A metamaterial system is called \emph{chiral} if its associated wallpaper group does not contain a reflection. A system is called \emph{achiral} otherwise.
    \end{definition}
    
    While this latter definition may not be worded the same as the description of the property found in the literature, it captures the same phenomenon (see Section~2.2 in \cite{review1}).

\subsection{Methods}
\label{sec:methods}
 
In order to faithfully compare auxetic behaviors across systems, we compare systems created by the same authors in the same project. This guarantees that we are comparing systems which used the same or similar parameters within their metamaterial systems, and all simulation was done with the same program and settings. Just as the Poisson's ratio of a system can be dependent on the directions of strain, the magnitude of strain can affect a system's auxeticity. While the study of deformation mechanisms of structures over large (finite) strains is important to applications of these systems, many of the experiments and applications involving auxetic systems focus on behavior over infinitesimal strains.
 
  \begin{convention}[Infinitesimal Strain and Finite Strain]
     We restrict our analysis to the initial (infinitesimal strain) Poisson's ratio of systems and denote it as $\nu_{xz}$ or $\nu_{zx}$. 
\end{convention}

\begin{figure}[H]
    \centering
  \includegraphics[width=.6\linewidth]{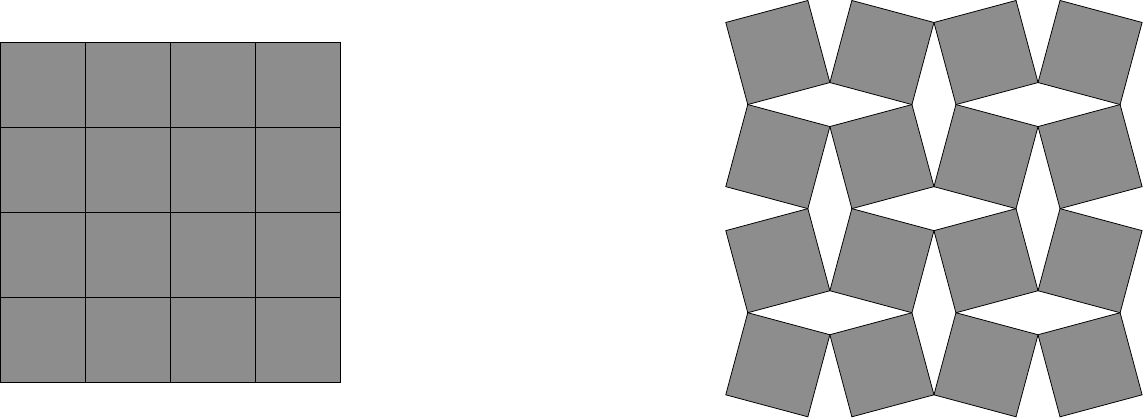}
\caption{The rotating squares system has different symmetries in its closed and partially open states}
\label{fig:partopen}
\end{figure}

    The symmetries of an auxetic system can potentially change with the state of the system. Consider, for example, the rotating squares model proposed by Grima and Evans \cite{grima2000auxetic} which is presented in Figure~\ref{fig:partopen}. In the completely closed state, the geometry is a regular tessellation of squares with wallpaper group $^*442$. Upon expansion, the system loses symmetry along its diagonals, changing its wallpaper group to $442$.
    
\begin{convention}[Wallpaper group associated to partially open state]
    We  define the wallpaper group of an auxetic system to be the one associated to the system's partially open state. 
\end{convention}

For a practical guide to determining the wallpaper group associated to a repeating pattern, refer to the flowchart on page 128 of \cite{flowchart}. For ease of reference, this flowchart is recreated in the appendix of this paper using orbifold notation instead of crystallographic notation.

    Some systems under consideration will be anisotropic and have two different Poisson's ratio listed by the source literature.
    
\begin{convention}[Anistropy vs. Isotropy]
 For anisotropic systems, the values $\nu_{xz}$, $\nu_{zx}$ will be \emph{averaged}, and the resulting value will be considered the Poisson's ratio of the system.
    \end{convention}

 \section{Data}
 \label{sec:data}
 
Our data come from two sources, \cite{korner} and \cite{comparative}. While it would clearly be preferable to have more examples of auxetic systems, the limitations of comparing like systems leave relatively few parts of the literature suitable for analysis. As mentioned in the previous section, in order to control for choices made by the authors, we will only compare systems coming from the same paper. In order to get meaningful analysis, we used sources that had examples with many different wallpaper groups. The papers listed here were the only ones we found that met these criteria, however this data is sufficient to suggest certain conclusions that we describe in Section~\ref{sec:conclusions}.

We briefly summarize the aims of each paper and the parameter choices that the authors made, before listing the data we extract from each.

 \subsection{Körner and Liebold-Ribeiro - \textit{A systematic approach to identify cellular auxetic metamaterials}}
 Körner and Liebold-Ribeiro \cite{korner} introduce a set of auxetic systems derived from the triangular, quadratic, and hexagonal lattices. In particular, they consider all the `vibrations', or sinusoidal perturbations, of these lattices and look at the special cases where all parts of the lattices move together at the same frequency (these special vibrations are called the \emph{eigenmodes} of the lattices). After specifying the parameters $L$ and $t$ below to be $L=2.5$mm for the hexagon, $L=5$mm for the square and triangle and $t=0.25$mm, the authors model the lattices in Abaqus 6.13 with the material properties of titanium and calculate the first eigenmodes of each of the lattices. 

 The authors proceed to create periodic lattices out of the eigenmodes and calculate the initial (infinitesmal) Poisson's ratio through simulations in Abaqus 6.13 for the resulting systems. Eight of these systems show auxetic behavior; these auxetic systems their wallpaper groups are given in the table below.
\vspace{0.2in}
 \begin{center}
\begin{tabular}{|| c c c ||} 
 \hline
   Picture of System & Wallpaper Group & Poisson's Ratio \cite{korner}\\ 
  \hline\hline
  \includegraphics[scale=0.3]{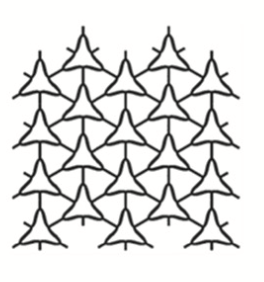}& $3^*3$ & $\nu_{xz}=\nu_{zx}=-0.2$\\[2ex] 
  \hline
   \includegraphics[scale=0.3]{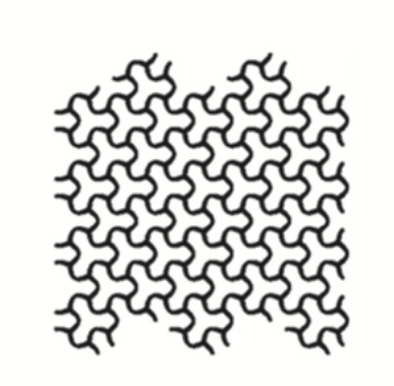}& $^*333$ & $\nu_{xz}=-0.1$, $\nu_{zx}=-0.3$\\[2ex]
  \hline
   \includegraphics[scale=0.3]{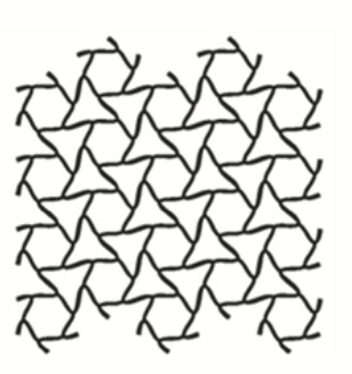}& $632$ & $\nu_{xz}=\nu_{zx}=-0.2$\\[2ex]
  \hline
   \includegraphics[scale=0.3]{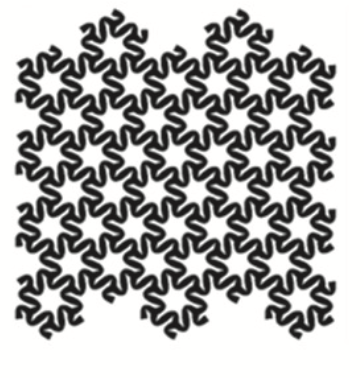}& $632$ & $\nu_{xz}=\nu_{zx}=-0.1$\\[2ex]
 \hline
   \includegraphics[scale=0.3]{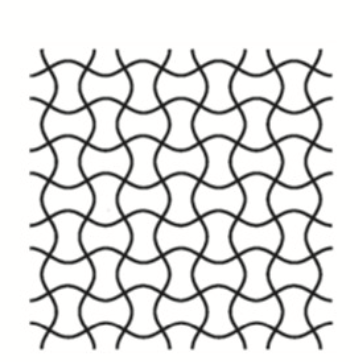}& $4^*2$ & $\nu_{xz}=\nu_{zx}=-0.3$\\[2ex]
  \hline
   \includegraphics[scale=0.3]{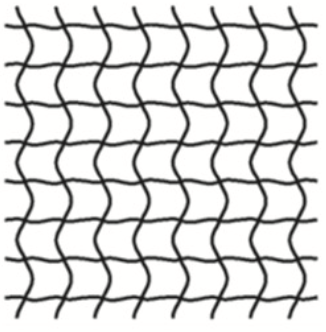}& $22^*$ & $\nu_{xz}=-1.6$, $\nu_{zx}=-0.2$\\[2ex]
  \hline
   \includegraphics[scale=0.3]{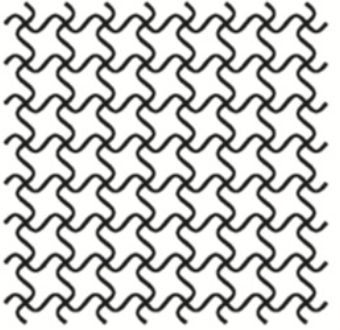}& $442$ & $\nu_{xz}=\nu_{zx}=-0.5$\\[2ex]
  \hline
   \includegraphics[scale=0.3]{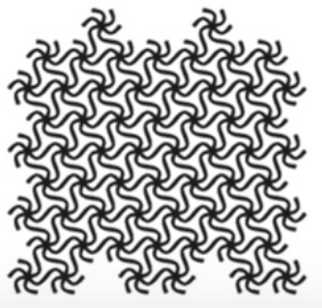}& $632$ & $\nu_{xz}=\nu_{zx}=-0.4$\\[2ex]
  \hline
 \end{tabular}
 \end{center}
 
 \newpage

 \subsection{Elipe and Lantada - \textit{Comparative study of auxetic geometries by means of computer-aided design and engineering} } 

Elipe and Lantada \cite{comparative} had the goal of creating a library of various auxetic systems. In order to do this, they consider a variety of different metamaterial systems and simulate their mechanical properties using the programs Solid Edge and NX-8.0 (Siemens PLM Solutions). To control for extraneous parameters, the authors attempt to make the dimensions of the systems as uniform as possible (see \cite[Section~2.1]{comparative}).
 
In the following table, we remark that the first system was considered in \cite{comparative} in two different configurations. While both configurations have the same wallpaper group, the version omitted from the table had a less optimal Poisson's ratio. We consider better version of this system only. Additionally, the authors also investigate a handful of systems with positive Poisson's ratio - we do not consider these systems in our data. In the table below we reproduce pictures of the systems considered by Elipe and Lantada.
 \vspace{0.2in}
  \begin{center}
 \begin{tabular}{|| c c c ||} 
 \hline
   Picture of System & Wallpaper Group & Poisson's Ratio \cite{comparative} \\ 
  \hline\hline

 \multirow{3}{5em}{\includegraphics[scale=0.08]{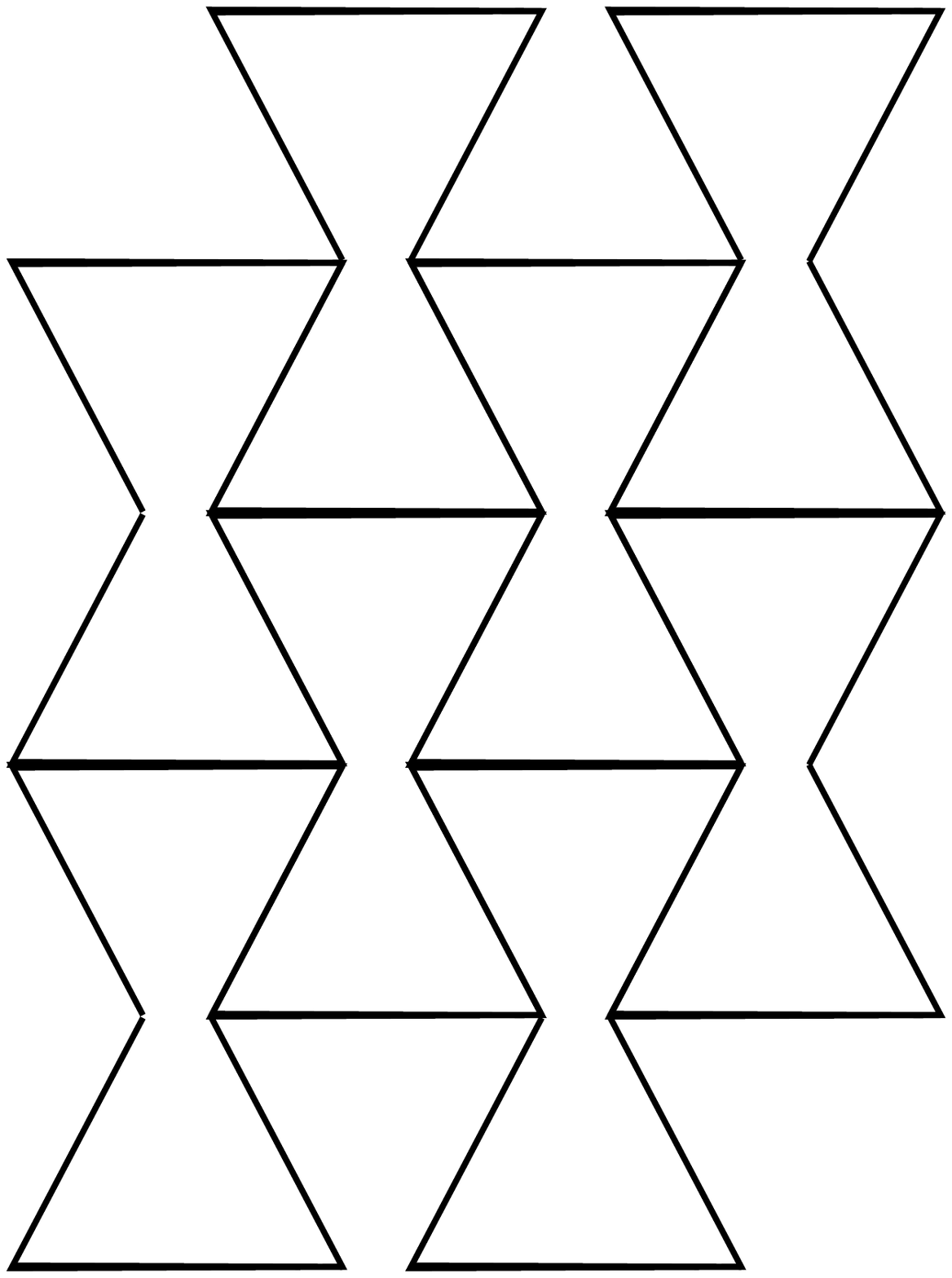}}& & \\ [2ex]
 & $2^*22$ &  $\nu_{xz}=-1.68$\\
 &  & \\
 
 \hline
 \multirow{3}{5em}{\includegraphics[scale=0.5]{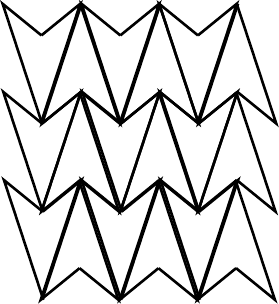}}& & \\ [2ex]
 & $22^*$ &  $\nu_{xz}=-0.789$\\
 &  & \\
 
  \hline
 \multirow{3}{5em}{\includegraphics[scale=0.1]{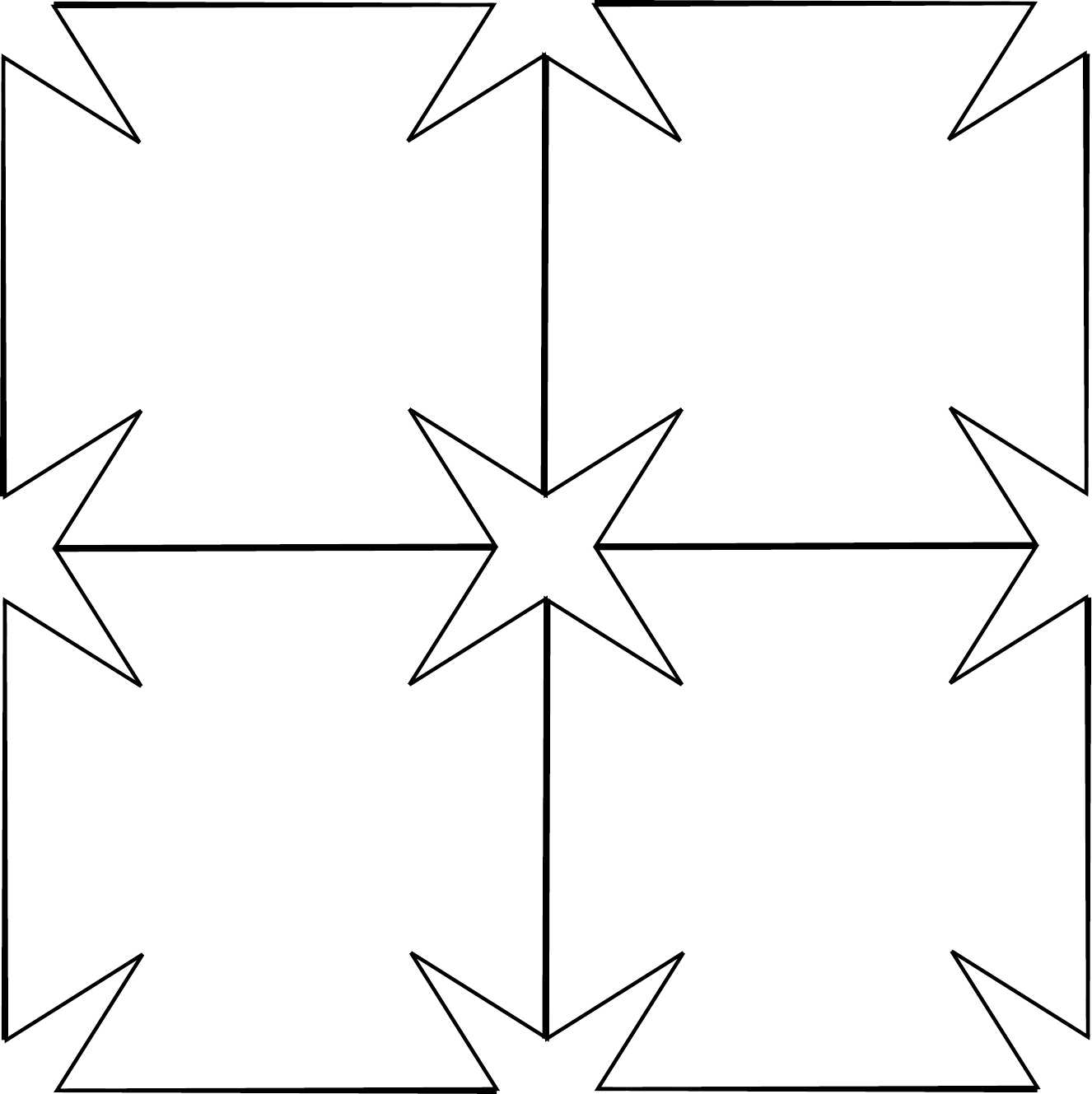}}& & \\ [2ex]
 & $^*442$ &  $\nu_{xz}=-0.504$\\
 &  & \\
   \hline
   
  \multirow{3}{5em}{\includegraphics[scale=0.1]{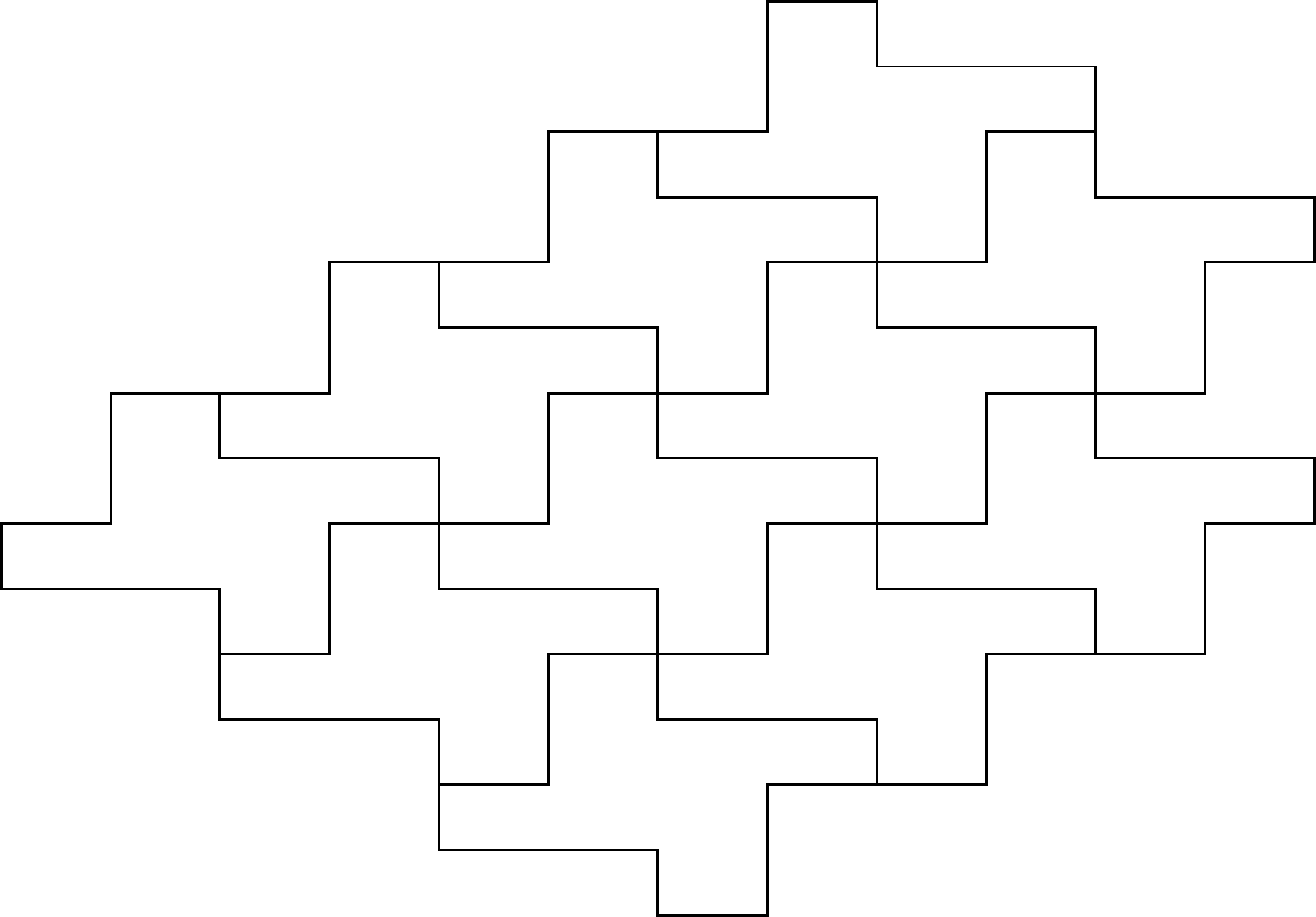}}& & \\ [2ex]
 & $2222$ &  $\nu_{xz}=-0.177$\\
 &  & \\
 \hline
 \multirow{3}{5em}{\includegraphics[scale=0.1]{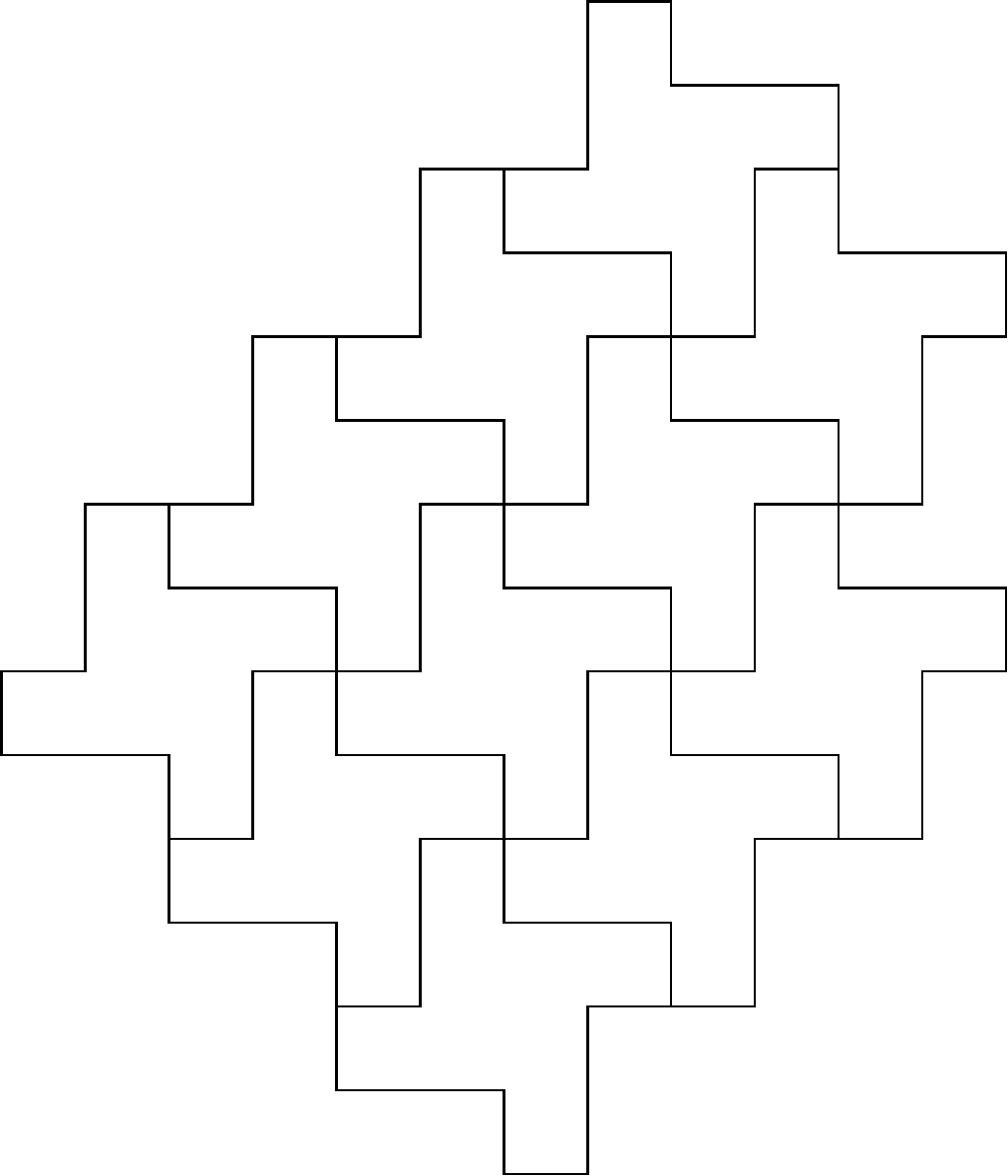}}& & \\ [2ex]
 & $442$ &  $\nu_{xz}=-0.326$ \\
  &  & \\
  \hline
  \multirow{3}{5em}{\includegraphics[scale=0.1]{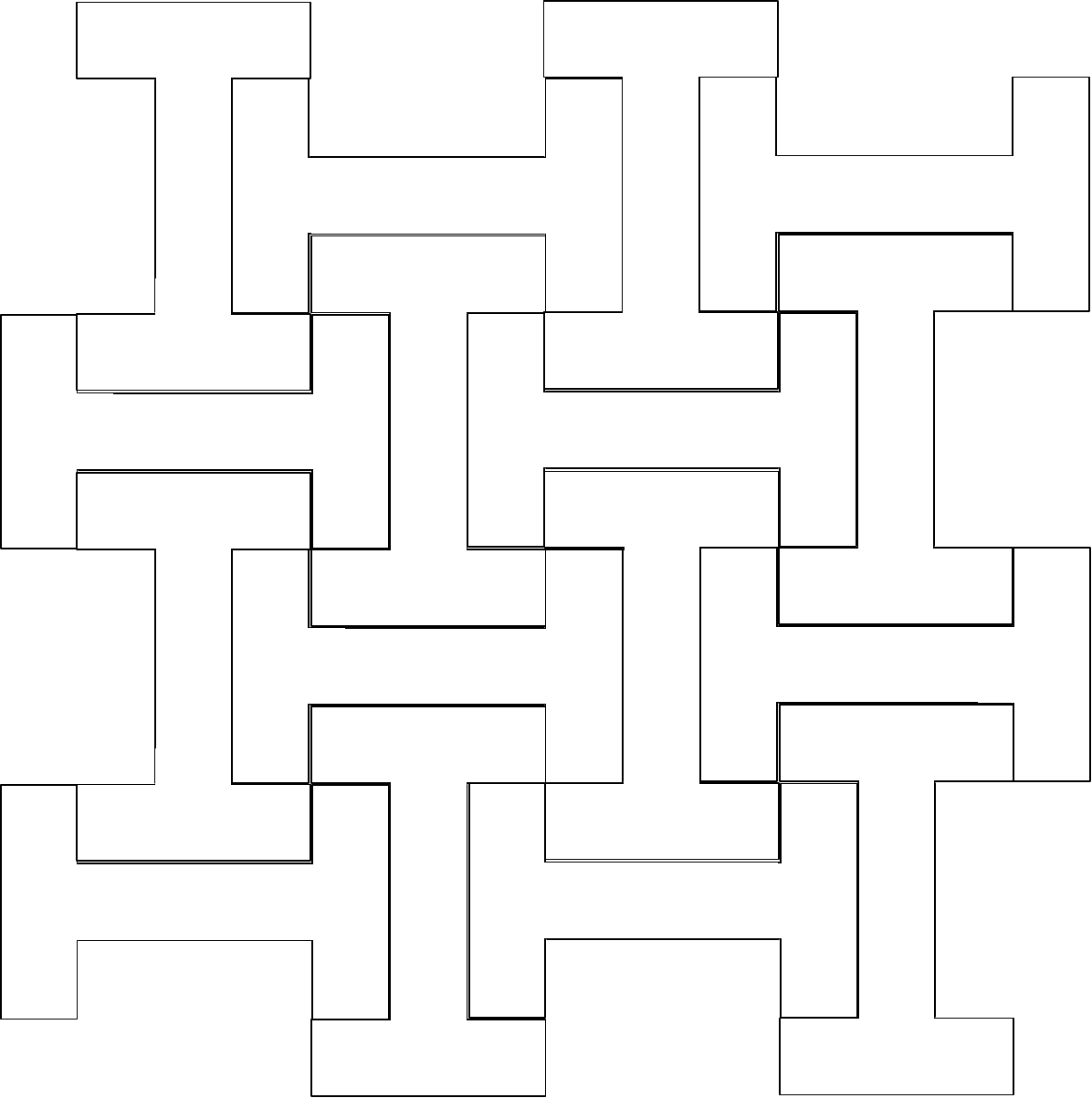}}& & \\  [2ex]
 & $4^*2$ & $\nu_{xz}=-0.901$\\
 &  & \\
 \hline
 \multirow{3}{5em}{\includegraphics[scale=0.07]{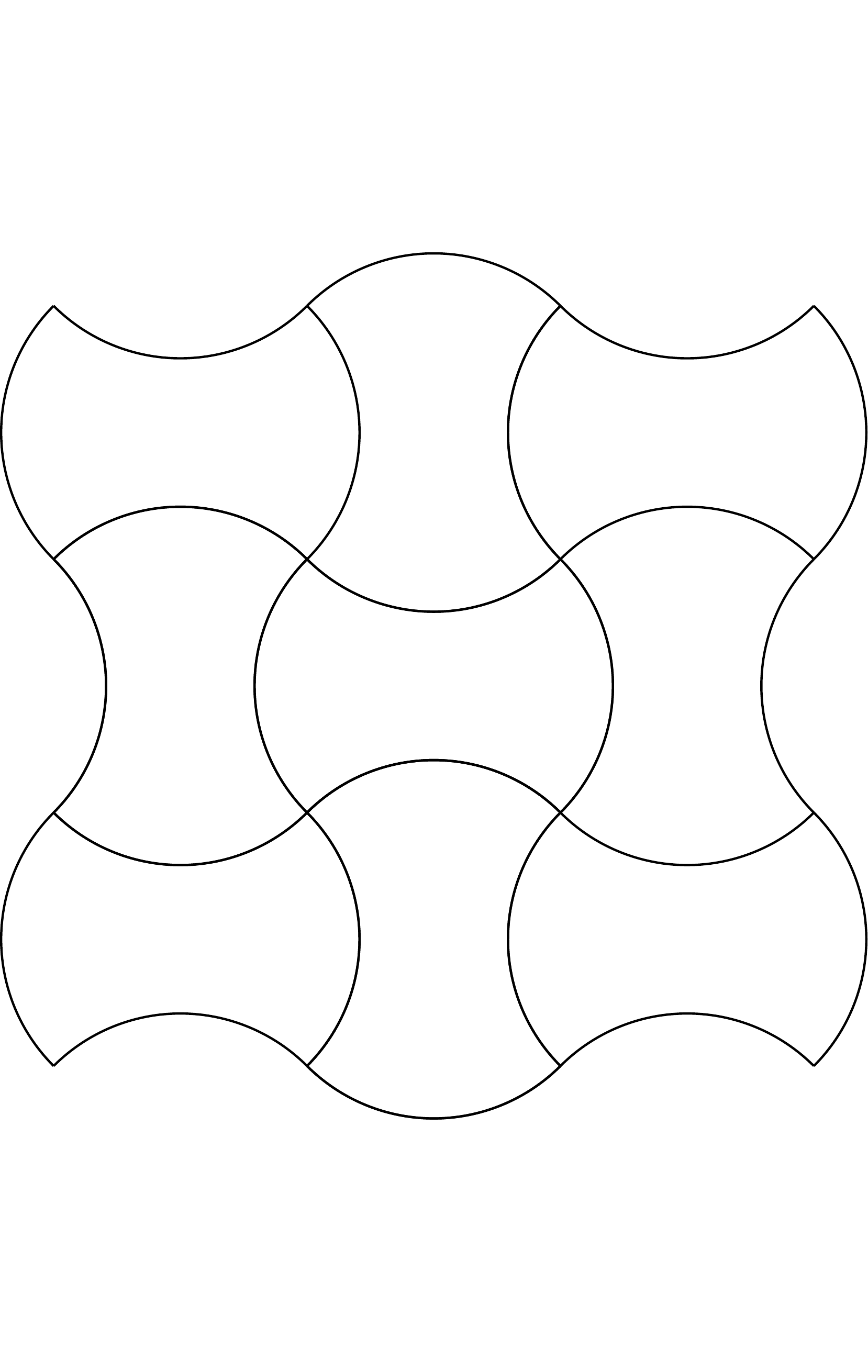}}& & \\  [2ex]
 & $4^*2$ & $\nu_{xz}=-0.81$\\
 &  & \\
 \hline
 \multirow{3}{5em}{\includegraphics[scale=0.09]{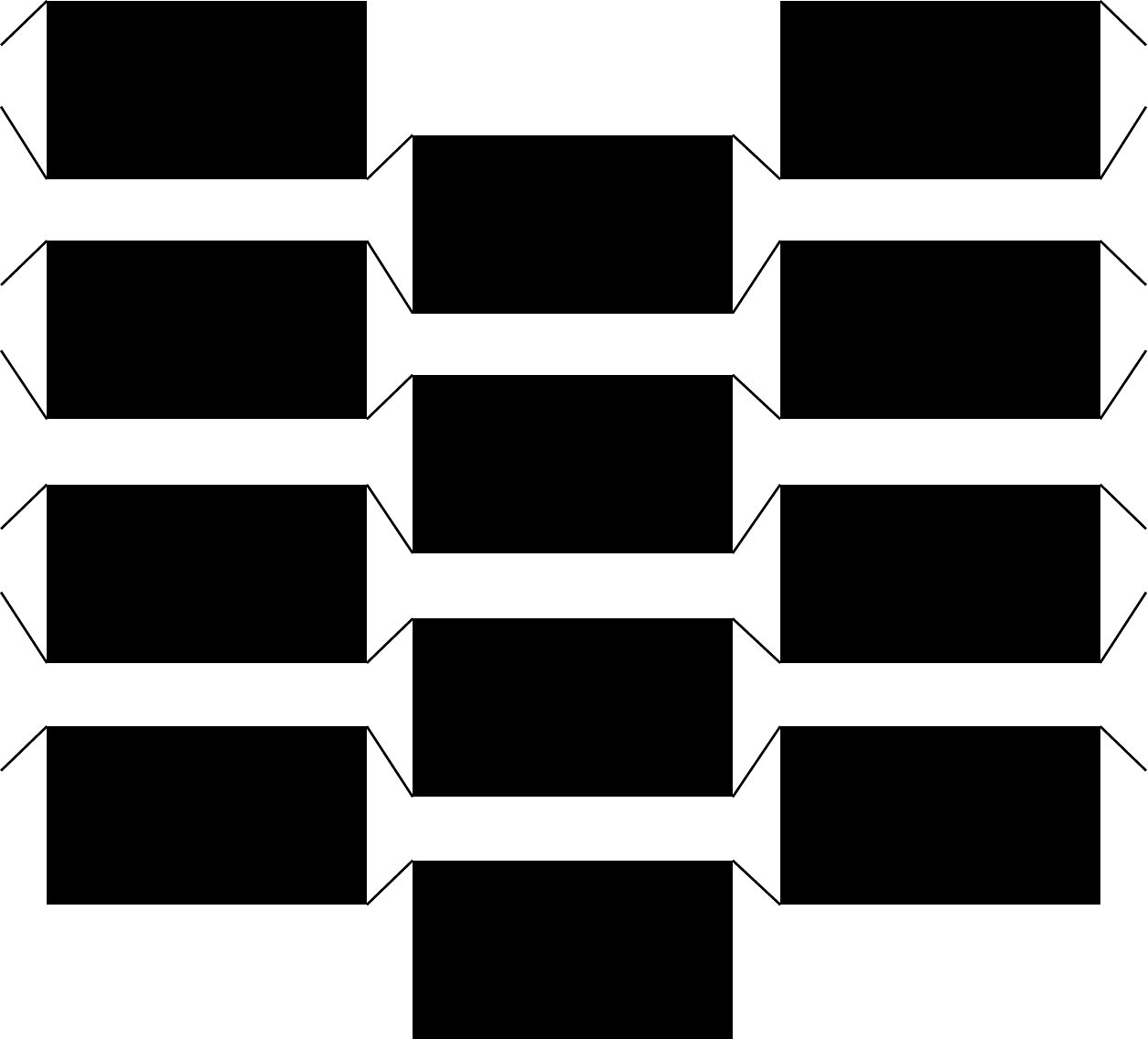}}& & \\  [2ex]
 & $2^*22$ & $\nu_{xz}=-0.071$ \\
 &  & \\
 \hline

\end{tabular}

\begin{tabular}{|| c c c ||} 
 \hline
   Picture of System & Wallpaper Group & Poisson's Ratio \cite{comparative}\\ 
  \hline 
  \hline
 \multirow{3}{5em}{\includegraphics[scale=0.12]{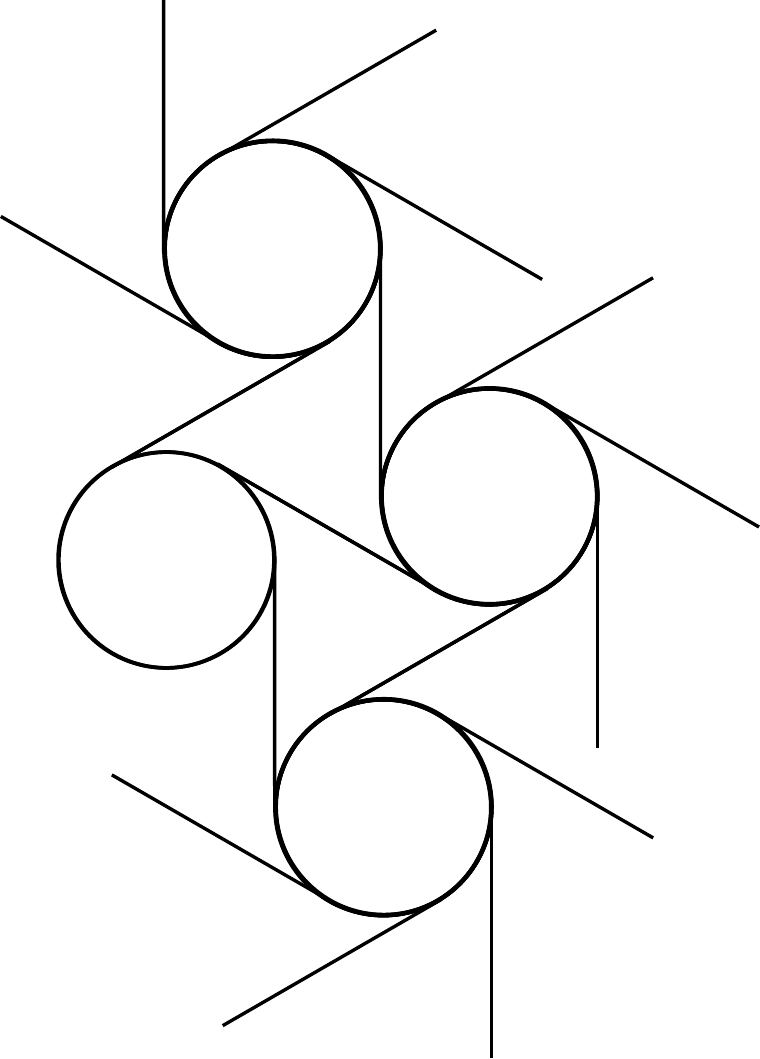}}& & \\  [2ex]
 & $632$ & $\nu_{xz}=-0.628$\\
 &  & \\
  \hline
\multirow{3}{5em}{\includegraphics[scale=0.12]{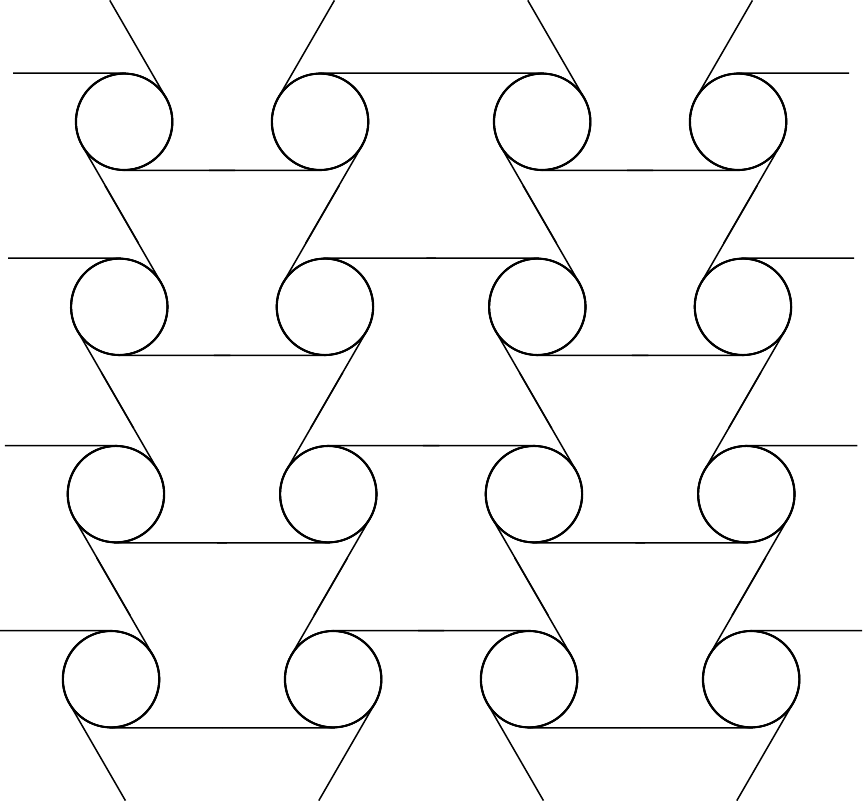}}& & \\  [2ex]
 & $22^*$ & $\nu_{xz}=-0.926$\\
 &  & \\

  \hline
\multirow{3}{5em}{\includegraphics[scale=0.1]{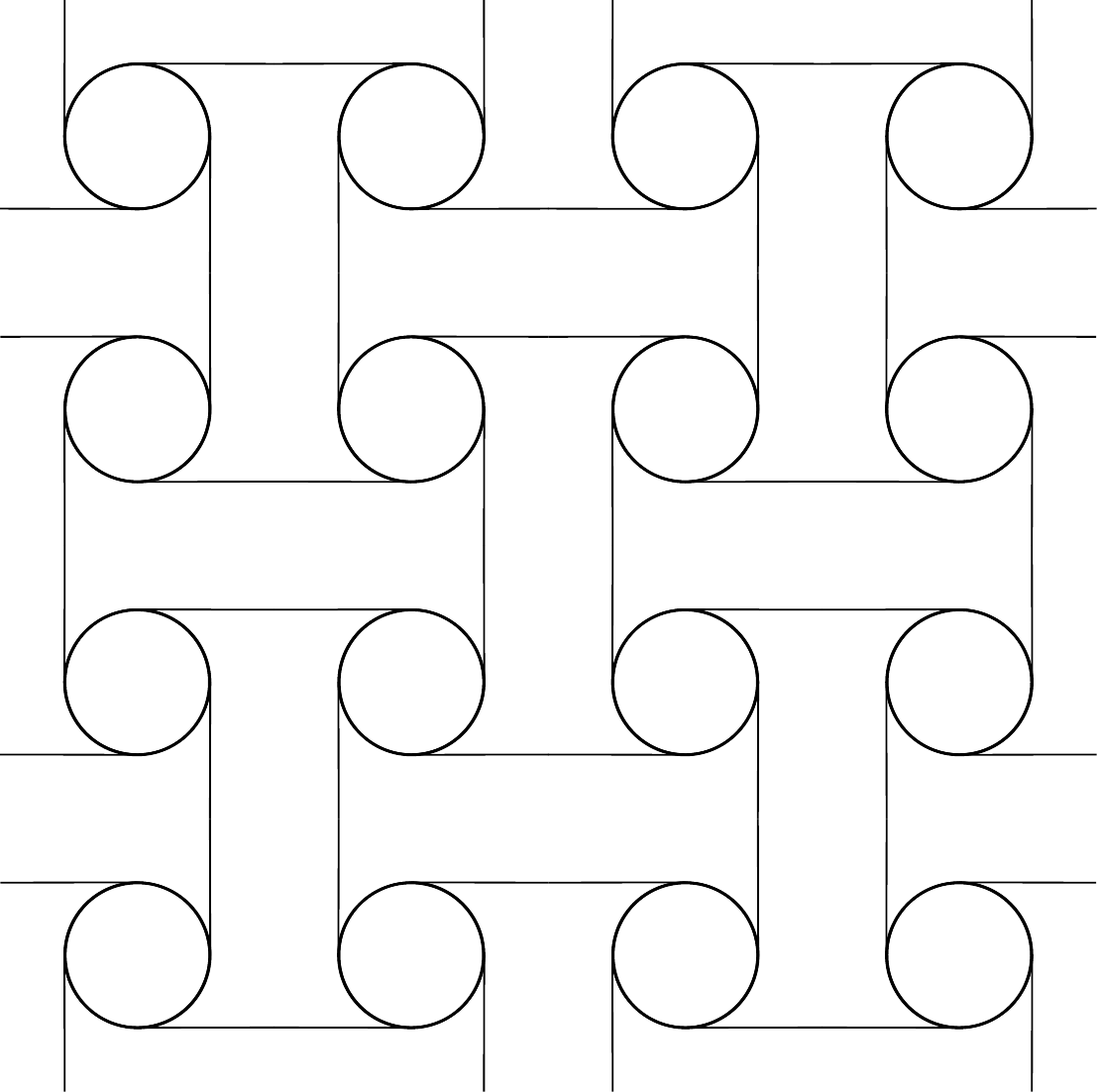}}& & \\  [2ex]
 & $4^*2$ & $\nu_{xz}=-0.966$\\
 &  & \\
 
  \hline

\multirow{3}{5em}{\includegraphics[scale=0.1]{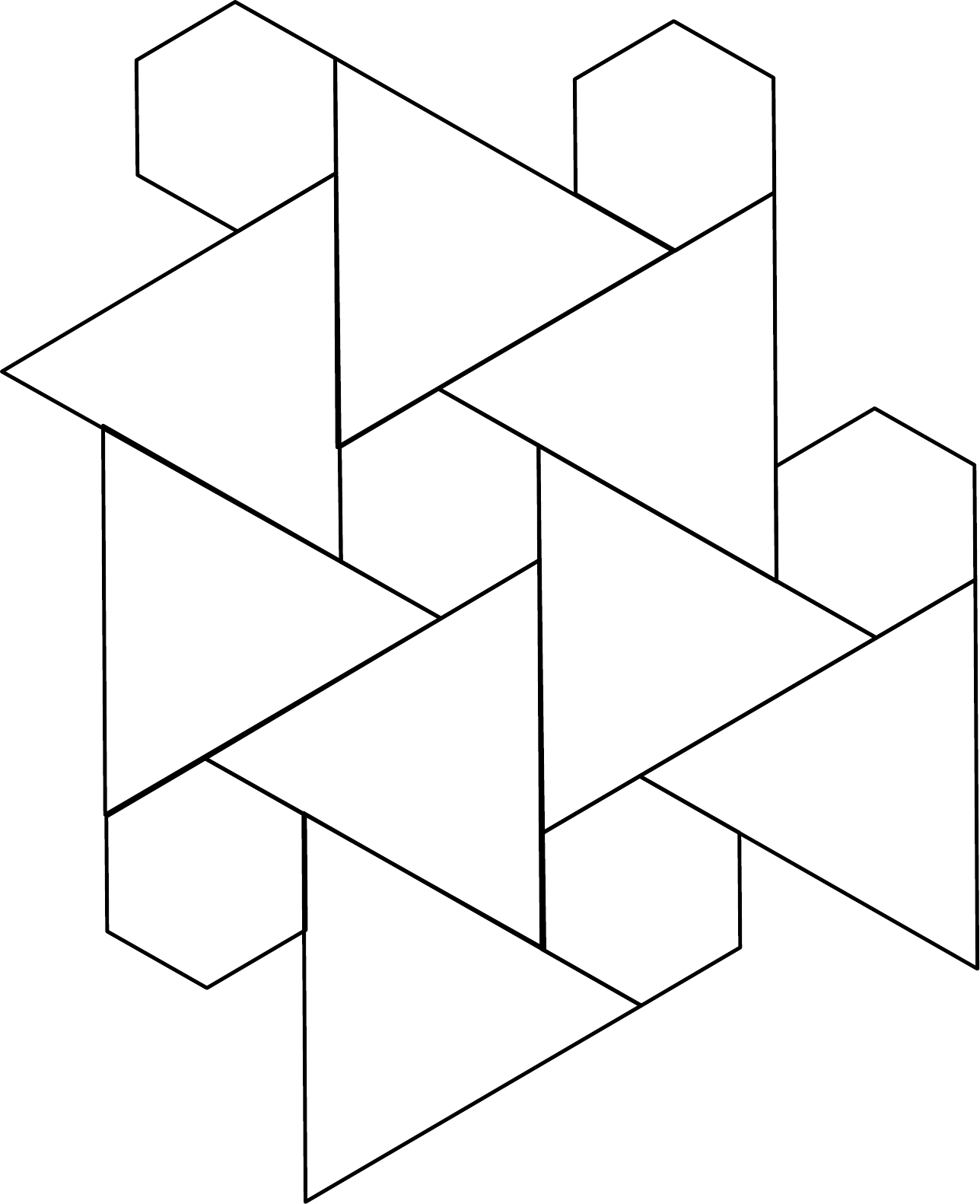}}& & \\  [2ex]
 & $632$ & $\nu_{xz}=-0.273$\\
 &  & \\
 \hline
\multirow{3}{5em}{\includegraphics[scale=0.1]{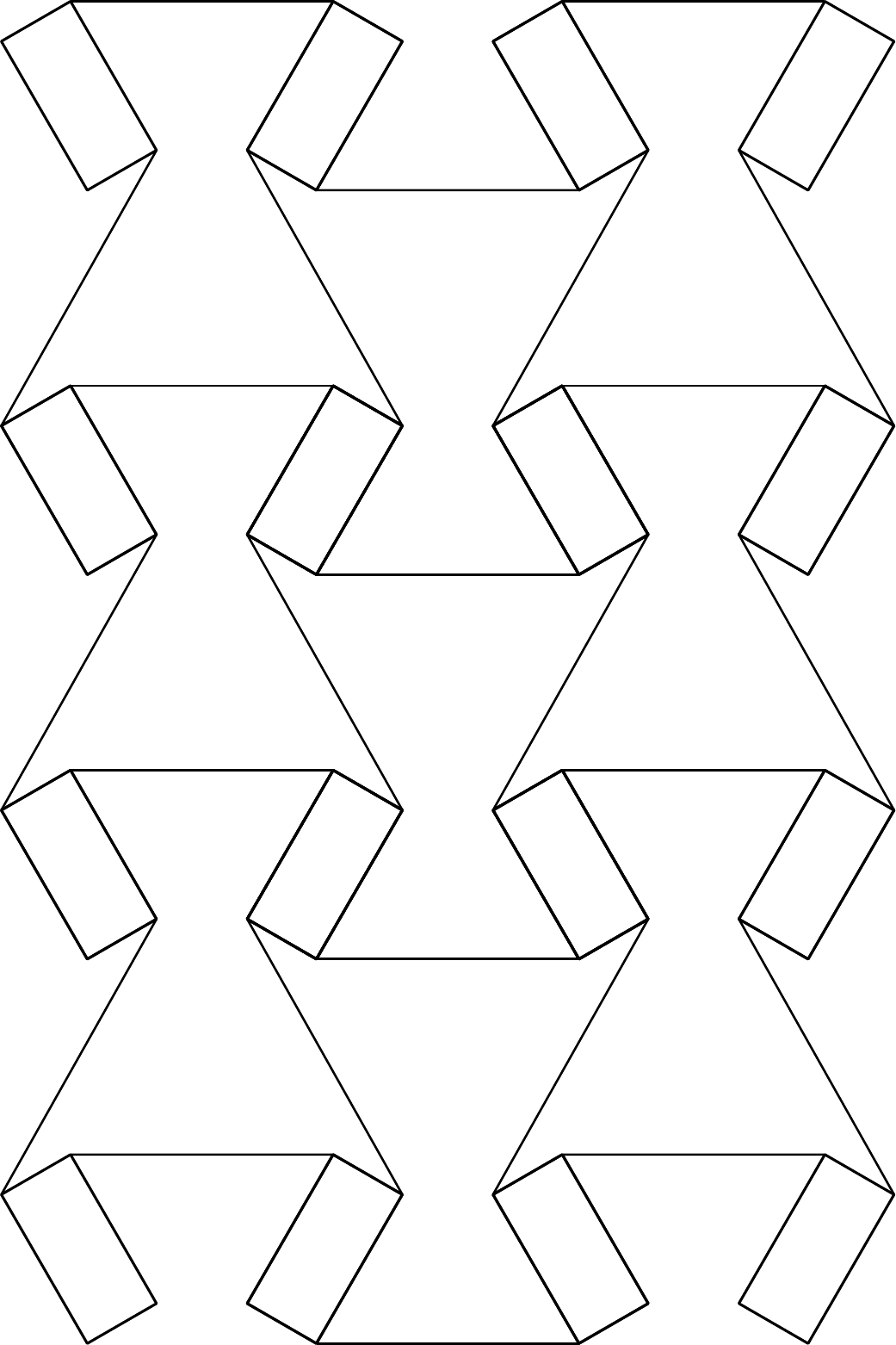}}& & \\  [2ex]
 & $22^*$ & $\nu_{xz}=-1.18$\\
 &  & \\  
 \hline
\multirow{3}{5em}{\includegraphics[scale=0.16]{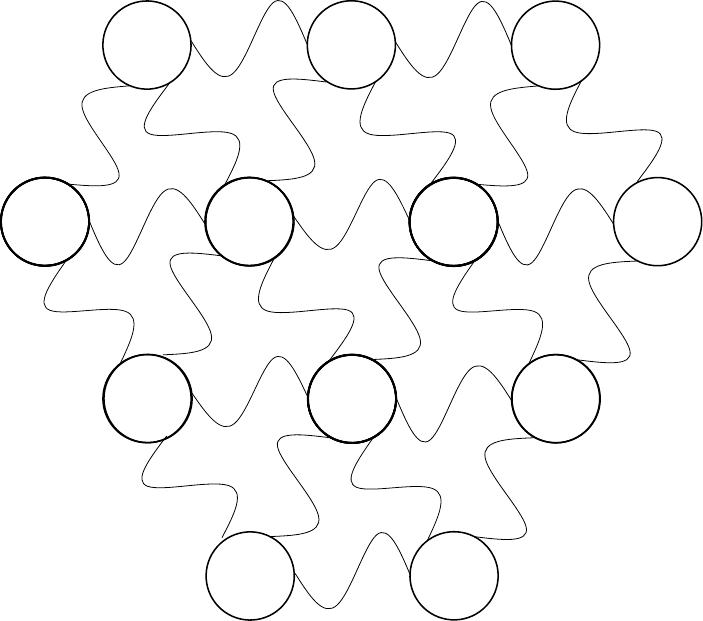}}& & \\  [2ex]
 & $632$ & $\nu_{xz}=-0.248$\\
 &  & \\
 \hline
\multirow{3}{5em}{\includegraphics[scale=0.22]{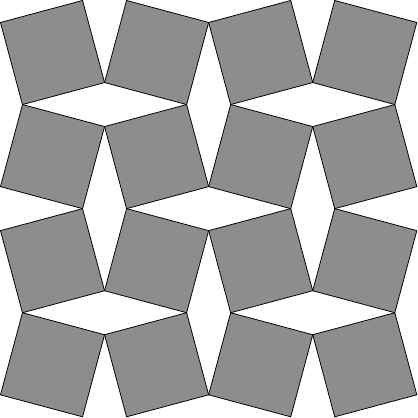}}& & \\  [2ex]
 & $4^*2$ & $\nu_{xz}=-0.315$\\
 &  & \\
 \hline
 \end{tabular}
  \end{center}

\section{Conclusions}
\label{sec:conclusions}
 We begin by listing the data from the previous section in a few suggestive ways. This first table compares the Poisson's ratio of systems which have reflections to those which do not.
 
 \vspace{0.2in}{}
 \begin{center}
 \begin{tabular}{|| l | l | l ||} 
 \hline
   Source & With reflection (achiral) & Without reflection (chiral)\\ 
  \hline\hline
  
  \cite{korner} & -0.2, -0.2, -0.3, -0.9 &  -0.2, -0.1, -0.5, -0.4\\
  \hline
    \cite{korner} averaged & \textbf{-0.4} & \textbf{-0.3} \\
    \hline
    \cite{comparative} &  -1.68, -0.789 -0.504 & -0.177, -0.326, -0.628\\ 
 & -0.901, -0.81, -0.071, -0.926 & -0.273, -0.248 \\ 
  & -0.966, -1.18, -0.315 & \\ 
    \hline
    \cite{comparative} averaged & \textbf{-0.8142} & \textbf{-0.330}\\
    \hline
  \end{tabular}
 \end{center}
 \vspace{0.2in}
 
 Based on this information, we posit that achiral auxetic systems exhibit lower Poisson's ratio than their chiral counterparts. This finding supports Hypothesis~\ref{hypothesis}.
  
  Next we consider the largest degree of rotation associated to each pattern. In examining the poset structure present in Figure~\ref{fig:poset}, we see that there are is a natural way to categorize the wallpaper groups by rotation. We will consider 6 and 3 fold rotations together, 4 fold rotations together, and everything that lies below both of these together. Because we did not observe any systems with no rotational symmetry whatsoever, this last class will include exactly those patterns with $180^\circ$ rotations.
  
  \vspace{0.2in}
   \begin{center}
   \begin{tabular}{|| l | l | l | l ||} 
 \hline
   Source & 6 \& 3 fold & 4 fold & 2 fold\\ 
  \hline\hline
  
  \cite{korner} & -0.2, -0.2, -0.2,  &  -0.3, -0.5 & -0.9\\
  
  & -0.1, -0.4 &   & \\
  \hline
    \cite{korner} averaged & \textbf{-0.22} & \textbf{-0.4} & \textbf{-0.9} \\
    \hline
    \cite{comparative} & -0.628,-0.273 &-0.504, -0.315 &  -1.68, -0.177 \\ 
 &-0.248 & -0.326, -0.901 & -0.789, -0.071
 \\
    & &  -0.81, -0.966 & -0.926,-1.18\\
    & & & \\ \hline
    \cite{comparative} averaged & \textbf{-0.383} & \textbf{-0.637} & \textbf{-0.804}\\
    \hline
  \end{tabular}
   \end{center}
  \vspace{0.2in}
  
  While the data in the first table supported Hypothesis~\ref{hypothesis}, this table tells a different story. In fact, the data suggests that \emph{fewer} rotations give better auxetic properties! This may be because there is more freedom to choose better parameters if one isn't restricted to a many-fold rotational structure. However, the authors conjecture that the `failure' of the higher fold patterns has more to do with the correlation with chirality. Indeed, in our data set, none of the 6-fold patterns were achiral (we found no examples of the wallpaper group *632), meaning that the collection of systems with 6-fold rotational symmetry were handicapped, in a certain sense.
  
  In \cite{comparative}, we can take a closer look at what is happening. Sorting these systems by their wallpaper group, we obtain the following:
  
  \vspace{0.2in}
   \begin{center}
   \begin{tabular}{|| l | l | l ||} 
 \hline
   Wallpaper Group & Poisson's Ratio & Averaged Poisson's Ratio\\ 
  \hline\hline
  $^*442$ & -0.504  &  \textbf{-0.504} \\
  \hline
   $4^*2$ & -0.81, -0.901, -0.966, -0.315 &  \textbf{-0.748}\\
  \hline
   $442$ & -0.326  & \textbf{-0.326}\\
  \hline
   $2222$ & -0.177 & \textbf{-0.177}\\
  \hline
   $2^*22$ & -1.68, -0.071  &  \textbf{-0.876} \\
  \hline
   $22^*$ & -0.789, -0.926, -1.18  & \textbf{-0.965}\\
  \hline
   $632$ & -0.273, -0.248, -0.628  & \textbf{-0.383} \\
  \hline
  \end{tabular}
   \end{center}
  \vspace{0.2in}
  
Looking at only the achiral systems and appealing to Figure~\ref{fig:poset}, we in fact see that less symmetric systems have lower (better) Poisson's ratios. But we can be even more specific. Systems with rotation centers not on reflection lines perform better than systems whose rotations are all on reflection lines. In fact, the wallpaper group with the lowest average Poisson's ratio had \emph{all} rotations not on reflection lines. We obtain a similar result from the systems in \cite{korner}; the achiral systems with some rotation centers not on reflection lines had average Poisson's ratio -0.467, while those having rotation centers all on reflection lines had an average Poisson's ratio of -0.2. This question is motivated by the discussion of rotational centers in \cite{korner} (see the observation regarding non-auxetic modes in Section~4.1).
  
We also note briefly that no systems in these papers had wallpaper groups appearing in the lower half (the bottom 5 nodes) of Figure~\ref{fig:poset}. While this lack of data does not allow us to conclude anything in particular about metamaterial systems with these wallpaper groups (if they exist), we believe that they are not found because they would not have more desirable properties than the systems currently found in the literature. The fact that the current literature on auxetic systems does not include any systems with these ``less symmetric'' wallpaper groups could be because they are not as desirable from a practical points of view. This suggests that more symmetric wallpaper groups are of more interest when it comes to producing auxetic metamaterials (up to the restriction described in the previous paragraph). 
  
Additionally, we found no wallpaper groups that had glide reflections (that were not compositions of translations and reflections that were themselves elements of the wallpaper group). Put another way, we never had to answer ``yes'' to the ``Glide reflection?'' question in the flowchart in Appendix A. Consequently, we may wonder if the groups $22 \times$, $^* \times$, and/or $\times \times$ can occur as wallpaper groups of auxetic systems.
  
  \begin{figure}
      \centering
  \includegraphics[scale=0.8]{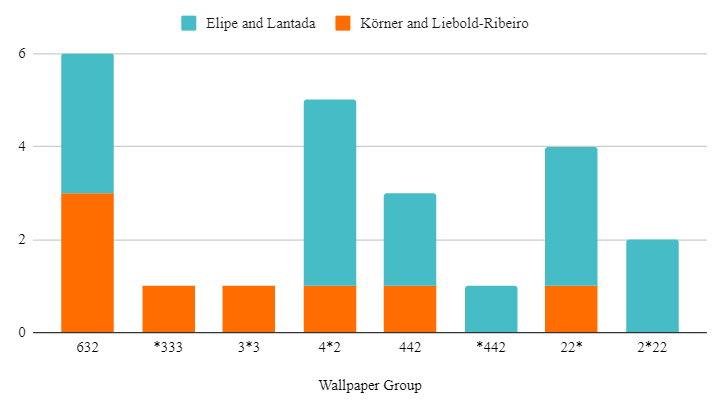}
      \caption{Frequency of wallpaper groups found in \cite{korner} and \cite{comparative}.}
      \label{fig:frequency}
  \end{figure}

  Figure~\ref{fig:frequency} shows that while the most common individual group in the two source papers has $6-$fold rotation, in aggregate $4-$fold patterns were more common.
 
 In summary, it appears that there are meaningful correlations between the mathematical structure (wallpaper group) of a two-dimensional auxetic system and the physical properties (Poisson's ratio) of that system. Our analysis suggests that an auxetic system whose wallpaper group contains at least one reflection is more likely to have a better (lower) Posisson's ratio than an auxetic system whose wallpaper group contains no reflections. It also suggests that auxetic systems whose wallpaper groups have $6$- or $3$-fold highest order rotations have worse (higher) Poisson's ratios than systems whose wallpaper groups have $4$- fold highest order rotations, which in turn have worse Poisson's ratios than systems whose wallpaper groups have $2$-fold highest order rotations. We acknowledge that our data are somewhat limited. However, these early results suggest that when designing auxetic metamaterial systems where a low Poisson's ratio is desirable, that it may be advantageous to focus on systems whose underlying geometric structure is a wallpaper group with reflections and $2$-fold rotations.

 \section{Further Questions}
 \label{sec:Further}
 
There are many further questions that could be asked about the role that geometry and topology play in the study of metamaterials. We list a few questions here that we believe would be particularly enlightening to investigate.

\begin{question}
  Some metamaterial design, such as re-entrant design (a design in which a repeating lattice of polygons or polyhedra is compressed so that all ribs protrude inward), is done using three-dimensional repeating patterns, rather than layering the same planar pattern at all heights (\cite{lakes1987foam}, \cite{3D1},\cite{3D2},\cite{3D3}). Could a similar analysis be done using space groups for these systems? For more information on space groups, the reader is referred to \cite{space_groups}.
\end{question}

\begin{question}
  Our analysis was done by comparing Poisson's ratio obtained using the infinitesimal strains for systems. Would similar patterns hold if one compared values associated to finite strains? We suspect the answer is yes. In this case, the amount of strain serves as an additional parameter to consider. It is possible that a more sophisticated mathematical relationship between these values could reveal deeper connections.
\end{question}

\begin{question}
Can the wallpaper group of a metamaterial system detect anisotropy? In the papers we analyzed, only two systems were anisotropic. These two systems had $2$- and $3$- fold maximum rotations. Does this pattern hold more broadly? For example, are systems with $2$-, $3$- or $6$- fold rotation more likely to be anisotropic than systems with $4$-fold rotation? 
\end{question}

\begin{question}
  Some mechanical engineering projects have obtained beneficial properties in systems by iterating metamaterial design at different scales, such as in \cite{hier}. This is closely related to the mathematical concept of a \emph{fractal}, a geometric object which is self-similar at various scales. Does the Hausdorff or Renyi dimension of the iteration of this scale-invariance make useful predictions about the behavior of an auxetic system?
\end{question}

\newpage

    \appendix

 \section{Wallpaper Group Resources}
 \label{sec:appendix}

\begin{adjustwidth}{}{}
\scalebox{0.85}{
 \noindent
\begin{tikzpicture}[
node distance=.8cm,
graynode/.style={rectangle, rounded corners, minimum height=1cm, text centered, draw=black!60, fill=black!5, very thick, minimum width=7mm},
greennode/.style={rectangle, rounded corners, minimum height=1cm, text centered, draw=green!60, fill=green!5, very thick, minimum width=7mm, text width=5em},
greennode1/.style={rectangle, rounded corners, minimum height=7mm, text centered, draw=green!60, fill=green!5, very thick, minimum width=7mm},
orangenode/.style={rectangle, rounded corners, minimum height=1cm, text centered, draw=orange!60, fill=orange!5, very thick, minimum width=7mm, text width=5em},
orangenode1/.style={rectangle, rounded corners, text centered, draw=orange!60, fill=orange!5, very thick, minimum width=7mm, minimum height=7mm,},
rednode/.style={rectangle, rounded corners, minimum height=1cm, text centered, draw=red!60, fill=red!5, very thick, minimum width=7mm,  text width=5em},
rednode1/.style={rectangle, rounded corners, minimum height=7mm, text centered, draw=red!60, fill=red!5, very thick, minimum width=7mm},
bluenode/.style={rectangle, rounded corners, minimum height=1cm, text centered, draw=blue!60, fill=blue!5, very thick, minimum width=7mm, text width=5em},
bluenode1/.style={rectangle, rounded corners, text centered, draw=blue!60, fill=blue!5, very thick, minimum width=7mm, minimum height=7mm,},
violetnode/.style={rectangle, rounded corners, minimum height=1cm, text centered, draw=violet!60, fill=violet!5, very thick, minimum width=7mm,  text width=5em},
violetnode1/.style={rectangle, rounded corners, text centered, draw=violet!60, fill=violet!5, very thick, minimum width=7mm, minimum height=7mm,},
]
\node[graynode]     (rotation)                                  {Highest order rotation?};
\node[]             (phantom1)        [below=of rotation]       {};

\node[bluenode]   (order3)        [below=of phantom1]         {Reflection?};
\node[]             (phantom2)      [below=of order3]           {};

\node[bluenode]   (order3rot)     [yshift=-1cm, xshift=1.5cm, below left=of order3]               {All rotation centers on reflection lines?};
\node[]             (phantom2a)      [below=of order3rot]           {};
\node[bluenode1]      (*333)          [below left=of phantom2a]       {$^*333$};
\node[bluenode1]      (3*3)           [below right=of phantom2a]      {$3^*3$};
\node[bluenode1]      (333)           [right=of 3*3]   {$333$};

\node[rednode]   (order4)        [left=of order3]                {Reflection?};
\node[]             (phantom3)      [below left=of *333]               {};

\node[rednode]   (order4rot)     [below left=of phantom3]                   {All rotation centers on reflection lines?};
\node[]             (phantom3a)      [below=of order4rot]           {};
\node[rednode1]      (*442)          [below left=of phantom3a]       {$^*442$};
\node[rednode1]      (4*2)           [below right=of phantom3a]      {$4^*2$};
\node[rednode1]      (442)           [right=of 4*2]             {$442$};

\node[orangenode]    (order6)    [left=of order4]                {Reflection?};
\node[]             (phantom4)      [below=of order6]           {};
\node[orangenode1]      (*632)      [below left= of phantom4]         {$^*632$};
\node[orangenode1]      (632)       [below right= of phantom4]        {$632$};

\node[violetnode]    (order2)    [right=of order3]               {Reflection?};
\node[]             (phantom5)      [xshift=2.5cm, below right=of 3*3]               {};
\node[violetnode]      (order2ref)           [below left=of phantom5]             {Reflection lines in two directions?};
\node[violetnode]   (order2glide)     [below right=of phantom5]                   {Glide reflection?};
\node[]             (phantom5a)      [below=of order2ref]           {};
\node[violetnode]      (order2refa)          [below left=of phantom5a]       {All rotation centers on reflection lines?};
\node[]             (phantom5b)      [below=of order2refa]           {};
\node[violetnode1]       (2*22)           [below right=of phantom5b]      {$2^*22$};
\node[violetnode1]       (*2222)           [left=of 2*22]      {$^*2222$};
\node[violetnode1]       (22*)           [right=of 2*22]      {$22^*$};
\node[violetnode1]       (22x)           [right=of 22*]      {$22\times$};
\node[violetnode1]      (2222)           [right=of 22x]      {$2222$};
\node[]             (phantom5c)      [xshift=2.65cm, right=of phantom5a]           {};

\node[greennode]    (order0)    [right=of order2]               {Reflection?};
\node[]             (phantom6)      [below=of order0]           {};
\node[greennode]      (order0g1)      [xshift=1cm, below left= of phantom6]         {Glide reflection?};
\node[greennode]       (order0g2)        [right= of order0g1]        {Glide reflection?};
\node[]             (phantom6a)      [below=of order0g1]           {};
\node[]             (phantom6b)      [below=of order0g2]           {};
\node[greennode1]       (*x)        [below left= of phantom6a]        {$^*\times$};
\node[greennode1]      (**)      [below right= of phantom6a]         {$^{**}$};
\node[greennode1]       (xx)        [right= of **]        {$\times\times$};
\node[greennode1]      (o)      [right= of xx]         {$\circ$};

\draw[-] (rotation.south) -- node[anchor=south] {} (phantom1.south);
\draw[->] (phantom1.south) -| node[yshift=-.42cm, anchor=east] {6} (order6.north);
\draw[->] (phantom1.south) -| node[yshift=-.42cm, anchor=east] {4} (order4.north);
\draw[->] (phantom1.south) -- node[anchor=east] {3} (order3.north);
\draw[->] (phantom1.south) -| node[yshift=-.42cm, anchor=east] {2} (order2.north);
\draw[->] (phantom1.south) -| node[yshift=-.42cm, anchor=east] {none} (order0.north);
\draw[-] (order3.south) -- (phantom2.south);
\draw[->] (phantom2.south) -| node[anchor=south] {no} (333.north);
\draw[->] (phantom2.south) -| node[anchor=south] {yes} (order3rot.north);
\draw[-] (order3rot.south) -- (phantom2a.south);
\draw[->] (phantom2a.south) -| node[anchor=south] {yes} (*333.north);
\draw[->] (phantom2a.south) -| node[anchor=south] {no} (3*3.north);
\draw[-] (order6.south) -- (phantom4.south);
\draw[->] (phantom4.south) -| node[anchor=south] {no} (632.north);
\draw[->] (phantom4.south) -| node[anchor=south] {yes} (*632.north);
\draw[-] (order4.south) -- (phantom3.south);
\draw[->] (phantom3.south) -| node[anchor=south] {no} (442.north);
\draw[->] (phantom3.south) -| node[anchor=south] {yes} (order4rot.north);
\draw[-] (order4rot.south) -- (phantom3a.south);
\draw[->] (phantom3a.south) -| node[anchor=south] {no} (4*2.north);
\draw[->] (phantom3a.south) -| node[anchor=south] {yes} (*442.north);
\draw[-] (order2.south) -- (phantom5.south);
\draw[->] (phantom5.south) -| node[anchor=south] {no} (order2glide.north);
\draw[->] (phantom5.south) -| node[anchor=south] {yes} (order2ref.north);
\draw[-] (order2ref.south) -- (phantom5a.south);
\draw[->] (phantom5a.south) -| node[anchor=south] {no} (22*.north);
\draw[->] (phantom5a.south) -| node[anchor=south] {yes} (order2refa.north);
\draw[-] (order2refa.south) -- (phantom5b.south);
\draw[->] (phantom5b.south) -| node[anchor=south] {no} (2*22.north);
\draw[->] (phantom5b.south) -| node[anchor=south] {yes} (*2222.north);
\draw[-] (order2glide.south) -- (phantom5c.south);
\draw[->] (phantom5c.south) -| node[anchor=south] {no} (2222.north);
\draw[->] (phantom5c.south) -| node[anchor=south] {yes} (22x.north);
\draw[-] (order0.south) -- (phantom6.south);
\draw[->] (phantom6.south) -| node[anchor=south] {no} (order0g2.north);
\draw[->] (phantom6.south) -| node[anchor=south] {yes} (order0g1.north);
\draw[-] (order0g1.south) -- (phantom6a.south);
\draw[->] (phantom6a.south) -| node[anchor=south] {no} (**.north);
\draw[->] (phantom6a.south) -| node[anchor=south] {yes} (*x.north);
\draw[-] (order0g2.south) -- (phantom6b.south);
\draw[->] (phantom6b.south) -| node[anchor=south] {no} (o.north);
\draw[->] (phantom6b.south) -| node[anchor=south] {yes} (xx.north);

\end{tikzpicture}}
\end{adjustwidth}

This chart is a guide for identifying the wallpaper group associated with any two-dimensional repeating pattern that fills the plane. It is based on a similar chart in \cite{flowchart}, but adapted to emphasize the orbifold structure and notation discussed in Section \ref{sec:Background}. The first question to ask is always ``What is the highest order rotation in the pattern?'' The order of a rotation is the minimum number of times the rotation must be applied before the pattern returns to the starting position. For example, if the smallest rotation is $90^\circ$, it is an order $4$ rotation because it takes four $90^\circ$ rotations to create a full $360^\circ$ rotation. The highest order rotation corresponds to the largest digit in the orbifold signature of the wallpaper group. The next question to ask is ``Are there any reflections?'' If there are one or more reflections, the orbifold signature of the wallpaper group will contain the $^*$ symbol. The question(s) to ask after determining the highest order rotation and whether or not there are any reflections depend on the previous answers. It should be noted that there are alternative questions that one could ask after the question about whether or not there are reflections.

\bibliographystyle{siamplain}
\bibliography{references.bib}

\begin{thebibliography}{10}

\bibitem{space_groups}
{\em International tables for crystallography. Vol. A: Space-group symmetry},
  Springer, 5. ed., reprinted with corrections~ed., 2005.

\bibitem{comparative}
{\sc J.~C. Alvarez~Elipe and A.~Diaz~Lantada}, {\em Comparative study of
  auxetic geometries by means of computer-aided design and engineering}, Smart
  Materials and Structures, 21 (2012), p.~105004,
  \url{https://doi.org/10.1088/0964-1726/21/10/105004}.

\bibitem{3D3}
{\sc S.~Babaee, J.~Shim, J.~C. Weaver, E.~R. Chen, N.~Patel, and K.~Bertoldi},
  {\em 3d soft metamaterials with negative poisson’s ratio}, Advanced
  Materials, 25 (2013), p.~5044–5049,
  \url{https://doi.org/10.1002/adma.201301986}.

\bibitem{borcea2015geometric}
{\sc C.~Borcea and I.~Streinu}, {\em Geometric auxetics}, Proceedings of the
  Royal Society A: Mathematical, Physical and Engineering Sciences, 471 (2015),
  p.~20150033.

\bibitem{streinu3}
{\sc C.~S. Borcea and I.~Streinu}, {\em Auxetic deformations and elliptic
  curves}, Computer Aided Geometric Design, 61 (2018), p.~9–19,
  \url{https://doi.org/10.1016/j.cagd.2018.02.003}.

\bibitem{streinu2}
{\sc C.~S. Borcea and I.~Streinu}, {\em Auxetic Regions in Large Deformations
  of Periodic Frameworks}, vol.~71, Springer International Publishing, 2019,
  p.~197–204, \url{https://doi.org/10.1007/978-3-030-16423-2_18},
  \url{http://link.springer.com/10.1007/978-3-030-16423-2_18}.

\bibitem{orbifold}
{\sc J.~H. Conway and D.~H. Huson}, {\em The orbifold notation for
  two-dimensional groups}, Structural Chemistry, 13 (2002), p.~247–257,
  \url{https://doi.org/10.1023/A:1015851621002}.

\bibitem{coxeterbook}
{\sc H.~S.~M. Coxeter}, {\em Introduction to geometry}, Wiley Classics Library,
  John Wiley \& Sons, Inc., New York, 1989.
\newblock Reprint of the 1969 edition.

\bibitem{generators}
{\sc H.~S.~M. Coxeter and W.~O.~J. Moser}, {\em Generators and Relations for
  Discrete Groups}, Springer Berlin Heidelberg, 1972,
  \url{https://doi.org/10.1007/978-3-662-21946-1},
  \url{http://link.springer.com/10.1007/978-3-662-21946-1}.

\bibitem{3D2}
{\sc L.~D’Alessandro, V.~Zega, R.~Ardito, and A.~Corigliano}, {\em 3d auxetic
  single material periodic structure with ultra-wide tunable bandgap},
  Scientific Reports, 8 (2018), p.~2262,
  \url{https://doi.org/10.1038/s41598-018-19963-1}.

\bibitem{Fedorov}
{\sc E.~Fedorov}, {\em Symmetry in the plane}, Proceedings of the Imperial St.
  Petersburg Mineralogical Society, 2 (1891), p.~345–390.
\newblock In Russian.

\bibitem{hier}
{\sc R.~Gatt, L.~Mizzi, J.~I. Azzopardi, K.~M. Azzopardi, D.~Attard, A.~Casha,
  J.~Briffa, and J.~N. Grima}, {\em Hierarchical auxetic mechanical
  metamaterials}, Scientific Reports, 5 (2015), p.~8395,
  \url{https://doi.org/10.1038/srep08395}.

\bibitem{grima2000auxetic}
{\sc J.~N. Grima and K.~E. Evans}, {\em Auxetic behavior from rotating
  squares},  (2000).

\bibitem{random}
{\sc J.~N. Grima, L.~Mizzi, K.~M. Azzopardi, and R.~Gatt}, {\em Auxetic
  perforated mechanical metamaterials with randomly oriented cuts}, Advanced
  Materials, 28 (2016), p.~385–389,
  \url{https://doi.org/10.1002/adma.201503653}.

\bibitem{review1}
{\sc P.~U. Kelkar, H.~S. Kim, K.-H. Cho, J.~Y. Kwak, C.-Y. Kang, and H.-C.
  Song}, {\em Cellular auxetic structures for mechanical metamaterials: A
  review}, Sensors, 20 (2020), p.~3132,
  \url{https://doi.org/10.3390/s20113132}.

\bibitem{korner}
{\sc C.~Körner and Y.~Liebold-Ribeiro}, {\em A systematic approach to identify
  cellular auxetic materials}, Smart Materials and Structures, 24 (2015),
  p.~025013, \url{https://doi.org/10.1088/0964-1726/24/2/025013}.

\bibitem{lakes1987foam}
{\sc R.~Lakes}, {\em Foam structures with a negative poisson's ratio}, Science,
  235 (1987), pp.~1038--1041.

\bibitem{lim}
{\sc T.-C. Lim}, {\em Auxetic materials and structures}, Springer, 2015.

\bibitem{review}
{\sc P.~Ma}, {\em A review on auxetic textile structures, their mechanism and
  properties}, Journal of Textile Science \& Fashion Technology, 2 (2019),
  \url{https://doi.org/10.33552/JTSFT.2019.02.000526},
  \url{https://irispublishers.com/jtsft/fulltext/a-review-on-auxetic-textile-structures-their-mechanism-and-properties.ID.000526.php}.

\bibitem{applications}
{\sc M.~Mir, M.~N. Ali, J.~Sami, and U.~Ansari}, {\em Review of mechanics and
  applications of auxetic structures}, Advances in Materials Science and
  Engineering, 2014 (2014), p.~1–17,
  \url{https://doi.org/10.1155/2014/753496}.

\bibitem{polya}
{\sc G.~Pólya}, {\em Xii. Über die analogie der kristallsymmetrie in der
  ebene}, Zeitschrift für Kristallographie - Crystalline Materials, 60 (1924),
  \url{https://doi.org/10.1524/zkri.1924.60.1.278},
  \url{http://www.degruyter.com/view/j/zkri.1924.60.issue-1-6/zkri.1924.60.1.278/zkri.1924.60.1.278.xml}.

\bibitem{3D1}
{\sc S.~M. Sajadi, C.~F. Woellner, P.~Ramesh, S.~L. Eichmann, Q.~Sun, P.~J.
  Boul, C.~J. Thaemlitz, M.~M. Rahman, R.~H. Baughman, D.~S. Galvão, and
  et~al.}, {\em 3d printing: 3d printed tubulanes as lightweight hypervelocity
  impact resistant structures (small 52/2019)}, Small, 15 (2019), p.~1970284,
  \url{https://doi.org/10.1002/smll.201970284}.

\bibitem{wallpaper_first}
{\sc M.~Stavric and A.~Wiltsche}, {\em Geometrical elaboration of auxetic
  structures}, Nexus Network Journal, 21 (2019), p.~79–90,
  \url{https://doi.org/10.1007/s00004-019-00428-5}.

\bibitem{flowchart}
{\sc D.~K. Washburn and D.~W. Crowe}, {\em Symmetries of culture: theory and
  practice of plane pattern analysis}, Univ. of Washington Press, 3. print~ed.,
  1998.

\end{thebibliography}

\end{document}